\documentclass[12pt,a4paper]{amsart}
\usepackage[utf8]{inputenc}
\usepackage[T1]{fontenc}
\usepackage{amsmath}
\usepackage{amsthm}
\usepackage{amsfonts}
\usepackage{amssymb}
\usepackage{mathrsfs}
\usepackage{url}
\usepackage{mathdots}
\usepackage[english]{babel}
\usepackage{geometry}
\usepackage{verbatim}
\usepackage{hyperref}
\usepackage{float}
\usepackage{enumerate}
\usepackage{tikz-cd}
\usepackage[nolabel]{showlabels}
\usepackage{extarrows}
\usepackage{comment}
\usepackage{bbm}

\hypersetup{
    bookmarks=true,         
    unicode=false,          
    pdftoolbar=true,        
    pdfmenubar=true,        
    pdffitwindow=false,     
    pdfstartview={FitH},    
		pdftitle={Motivic weight {$23$}},    
    pdfauthor={Ga\"etan Chenevier and Olivier Ta\"ibi},     
    colorlinks=true,       
    linkcolor=blue,          
    citecolor=green,        
    filecolor=green,      
    urlcolor=cyan}           

\newcommand{\Q}{\mathbb{Q}}

\newcommand{\R}{\mathbb{R}}
\providecommand{\C}{\mathbb{C}}
\renewcommand{\C}{\mathbb{C}}

\newcommand{\Z}{\mathbb{Z}}

\newcommand{\ps}{\par \smallskip}
\newcommand{\isomo}{\overset{\sim}{\rightarrow}}

\newtheorem{thm}{Theorem}[section]

\newtheorem{lemma}[thm]{Lemma}

\newtheorem{coro}[thm]{Corollary}
\newtheorem{definition}[thm]{Definition}
\newtheorem{propdefi}[thm]{Definition-Proposition}
\newtheorem{prop}[thm]{Proposition}

\newtheorem{remark}[thm]{Remark}
\newtheorem{example}[thm]{Example}

\newenvironment{pf}
{\medskip\noindent {\it Proof.  }}
{\hfill\nobreak $\Box$ \par\bigbreak}

\begin{document}

\author{Ga\"etan Chenevier}

\author{Wee Teck Gan}

\address{Ga\"etan Chenevier, CNRS, D.M.A., \'Ecole Normale Sup\'erieure, 45 rue d'Ulm, 75005 Paris, France}
\email{gaetan.chenevier@math.cnrs.fr}

\address{Wee Teck Gan, Department of Mathematics, National University of Singapore, 10 Lower Kent Ridge Road
Singapore 119076} 
\email{matgwt@nus.edu.sg}

\thanks{The first author is supported by the C.N.R.S. and by the project ANR-19-CE40-0015-02 (COLOSS). 
The second author is partially supported by a Singapore government MOE Tier 1 grant R-146-000-320-114.
The authors thank Michael Larsen and Jun Yu for some useful discussions concerning their articles \cite{larsen1} and \cite{yu}.
}

\begin{abstract} 
We classify the pairs of group morphisms $\Gamma \rightarrow {\rm Spin}(7)$ which are element conjugate but not globally conjugate.  As an application, we study the case where $\Gamma$ is the Weil group of a
$p$-adic local field, which is relevant to the recent approach 
to the local Langlands correspondence for ${\rm G}_2$ and ${\rm PGSp}_6$ 
in \cite{gansavin}. As a second application, we improve some result in \cite{kretshin} about ${\rm GSpin}_7$-valued Galois representations.
\end{abstract}

\title{${\rm Spin}(7)$ is unacceptable}

\maketitle

\section{Introduction}

Let $G$ be a compact group and let $\Gamma$ be an arbitrary group. 
A group morphism $r : \Gamma \rightarrow G$ is called {\it unacceptable} 
if there is a morphism $r' :  \Gamma \rightarrow G$ such that: \ps\ps
(U1) for all $\gamma \in \Gamma$, $r(\gamma)$ and $r'(\gamma)$ are conjugate in $G$, \ps\ps
(U2) $r'$ is not $G$-conjugate to $r$, \ps\ps
\ps\ps
\noindent 
Otherwise, we say that $r$ is {\it acceptable}. 
If $r$ is unacceptable, then so is any $G$-conjugate of $r$. 
Following Larsen \cite{larsen1} we also say that $G$ is {\it acceptable} 
if every $G$-valued group morphism is acceptable. 
If $r$ and $r'$ satisfy condition (U1) we also say that they are {\it element conjugate} (in $G$). 
Since the characters of (finite dimensional, continuous) representations of $G$ separate its conjugacy classes,
condition (U1) is also equivalent to:\ps\ps

(U1)' for all representations 
$\rho : G \rightarrow {\rm GL}_n(\C)$, the representations $\rho \circ r$ and 
$\rho \circ r'$ of $\Gamma$ are isomorphic.\ps\ps

Beyond the group theoretic legitimacy of these notions, 
a motivation for their study comes from the theory of automorphic representations.
Indeed, it was observed long ago by Langlands that,
for a given reductive group $H$ over a number field, 
the unacceptability of (the compact form of) the $L$-group of $H$ 
creates serious local-global difficulties in the study of automorphic representations of $H$,
and in particular non multiplicity one phenomena (see  {\it e.g.} 
\cite{blasius} for the first instance of such a phenomenon). 
In a similar vein, acceptability questions typically arise when one tries to characterize 
a representation of the absolute Galois group of $\Q$ with values in a reductive group 
by its Frobenius conjugacy classes.
We also mention that 
these notions have applications to constructions 
of isospectral, non isometric, Riemannian manifolds (see the introduction in \cite{larsen1}).
\ps\ps 

It is a folklore result that for each integer $n\geq 1$, the classical compact groups ${\rm U}(n)$, 
as well as ${\rm O}(n)$ and ${\rm Sp}(n)$, are all acceptable \cite{larsen1}. 
It follows that ${\rm SU}(n)$ and that ${\rm SO}(2n+1)$ are acceptable as well. 
Also, it is not difficult to show that ${\rm SO}(4)$ is acceptable, 
and results of Griess in \cite{griess} also imply that ${\rm G}_2$ is acceptable. 
It follows from these facts and standard exceptional isomorphisms 
that ${\rm Spin}(n)$ is acceptable for $n \leq 6$. 
Our first result is that this fails for $n=7$. 

\begin{prop} 
The group ${\rm Spin}(n)$ is not acceptable for $n\geq 7$
\end{prop}

A classification of the compact connected acceptable Lie groups was 
initiated by Larsen in \cite{larsen1} and \cite{larsen2}, and more recently 
pursued by J. Yu \cite{yu}. For the little story, both \cite{larsen1} and the 
first version of \cite{yu} contained 
two different incorrect ``proofs'' that ${\rm Spin}(7)$ is acceptable! 
Some problem in Larsen's proof was discovered by the second author in 2017, 
after he wrote a manuscript on a local Langlands correspondence
for ${\rm PGSp}_6$ over $p$-adic fields, and in which the acceptability of ${\rm Spin}(7)$ 
played a crucial role.\footnote{The content of that manuscript just appeared as the Appendix C 
of the preprint~\cite{gansavin}.} A counterexample was then found by the first author.
Our aim in this paper is not only to explain this counterexample, 
by partially reproducing a letter we sent some time ago to Larsen and Yu \cite{letter}, but also 
to give a classification of the unacceptable ${\rm Spin}(7)$-valued morphisms,
with in view some applications to the local Langlands correspondence
for ${\rm G}_2$ and ${\rm PGSp}_6$ \cite{gansavin}. We end by stating a remarkable result proved by Yu in \cite{yu}:
a compact connected Lie group is acceptable if, and only if, its derived subgroup is a direct product
of the aforementioned acceptable groups. Also, we mention a study in \cite{wang} of the unacceptable
continuous morphisms $\Gamma \rightarrow {\rm SO}(2n)$ with $\Gamma$ compact and connected. \ps

A first general result we prove in Section~\ref{sect:generaln}, which holds for all odd $n$, is the following.
We denote by $E$ the standard representation of ${\rm Spin}(n)$, 
an $n$-dimensional {\it real} representation,
and by $S$ a Spin representation, a $2^{k}$-dimensional {\it complex} representation with $n=2k+1$.
If $r : \Gamma \rightarrow {\rm Spin}(n)$ is a given morphism, 
it allows to view $E$ as an $\R[\Gamma]$-module,
and $S$ as a $\C[\Gamma]$-module. \ps

\begin{thm} 
\label{thm:gen}
A morphism $r : \Gamma \rightarrow {\rm Spin}(n)$, with $n$ odd, is unacceptable if, and only if,
there is an order $2$ character $\eta : \Gamma \rightarrow \{ \pm 1\}$ such that: 
\begin{itemize}
\item[(i)] we have $S \simeq S \otimes \eta$ as $\C[\Gamma]$-modules,\ps
\item[(ii)] no $\R[\Gamma]$-submodule of $E$ has determinant $\eta$.
\end{itemize}
\end{thm}

Of course, if $\Gamma$ has no order two character, 
then any morphism $\Gamma \rightarrow {\rm Spin}(n)$ is acceptable (an easier fact).
Although the theorem above is interesting, it does not explain the shape of the unacceptable
morphisms (nor why they only exist for $n\geq 7$).
Our second type of results are thus more specific to the case $n=7$.
A subgroup of ${\rm Spin}(7)$ is called a {\it ${\rm Spin}(1,6)$-subgroup} 
if it is obtained as the inverse image, via the canonical map 
${\rm Spin}(7) \rightarrow {\rm SO}(7)$, of the stabilizer of a line in $\R^7$. 
All these subgroups are conjugate in ${\rm Spin}(7)$ and are semi-direct 
products of $\Z/2\Z$ by ${\rm Spin}(6) \simeq {\rm SU}(4)$; however,
they are not isomorphic to ${\rm Pin}(6)$. The main result of Section~\ref{sect:n7} is:
\ps\ps


\begin{thm} 
\label{thm:exceptionPin6}
Assume $r : \Gamma \rightarrow {\rm Spin}(7)$ is unacceptable. 
Then its image $r(\Gamma)$ is contained in a ${\rm Spin}(1,6)$-subgroup of ${\rm Spin}(7)$,
or equivalently, the $\R[\Gamma]$-module $E$ contains a stable line.
\end{thm}

We stress that the necessary condition above is by no mean sufficient for being unacceptable: 
the inclusion of a ${\rm Spin}(1,6)$-subgroup in ${\rm Spin}(7)$ is acceptable.
Nevertheless, if $r : \Gamma \rightarrow {\rm Spin}(7)$ is unacceptable then $E$ contains some character 
$\chi : \Gamma \rightarrow \{ \pm 1\}$. We say that $r$ is {\it of type} I if we may take $\chi =1$.
Otherwise, we let $\Gamma_0 \subset \Gamma$ be the kernel of $\chi$ and say that 
$r$ is {\it of type} II (with respect to $\chi$) if $r_{|\Gamma_0}$ is unacceptable, and {\it of type} III otherwise (the precise definitions are slightly more constraining: see Definitions~\ref{def:deftypeII} and~\ref{def:typeIII}.) 
The main goal of the remaining parts of the paper is to give a classification of the unacceptable morphisms of each type. Essentially, we first give in each case some examples and then 
show that they are universal in some precise sense. 
Three specific compact subgroups of ${\rm Spin}(7)$, that we denote by 
${\mathcal{G}}$, $\mathcal{H}$ and $\mathcal{I}$ in the paper, and whose given
embeddings into ${\rm Spin}(7)$ are all unacceptable, play an important role. 
They are extensions of $\Z/2\Z \times \Z/2\Z$ by ${\rm SO}(2) \times {\rm SO}(2)$, ${\rm SO}(2) \times {\rm SO}(2) \times {\rm SO}(2)$ and ${\rm SO}(4)$  respectively. \ps

In Section~\ref{sect:typeI}, we prove that up to conjugacy, 
any unacceptable morphism of type {\rm I} 
factors through ${\mathcal{G}}$ (Theorem~\ref{thm:univGG}). 
Conversely, we give two necessary and sufficient conditions 
for a morphism into ${\mathcal{G}}$ to give rise to 
an unacceptable morphism into ${\rm Spin}(7)$ 
(Propositions~\ref{prop:converseN} and~\ref{prop:classnonex}). 
In Section~\ref{sect:typeII}, we reduce the study 
of type {\rm II} morphisms to that of type {\rm I} ones. 
Then we prove in Section~\ref{sect:typeIII} that up to conjugacy, 
any unacceptable morphism of type {\rm III} factors through 
either $\mathcal{H}$ or $\mathcal{I}$ 
(see Theorems~\ref{thm:exclasstypIIIa} and~\ref{thm:exclasstypIIIb}). 
The two corresponding situations are called 
type {\rm IIIa} and {\rm IIIb} respectively. A few lemmas needed in the proofs 
are gathered in the Appendix. A simple corollary of these results is:

\begin{coro}
Assume $r : \Gamma \rightarrow {\rm Spin}(7)$ is unacceptable.
Then up to conjugating $r$ if necessary, 
the image of $r$ is either included in ${\mathcal{G}}$, $
\mathcal{H}$ or $\mathcal{I}$, or has an index 2 subgroup included in ${\mathcal{G}}$.
\end{coro}

We leave as an open problem the question of classifying 
the unacceptable ${\rm Spin}(n)$-valued morphisms for $n>7$. \ps
In Section~\ref{sec:weilcase} we finally study the special case where $\Gamma$ is the Weil group 
${\rm W}_F$ of a finite extension $F$ of $\Q_p$. The general question, of inverse Galois theory flavour, 
is to understand to what extend the variety of general examples discussed above 
does occur for ${\rm W}_F$. 
For instance, as we shall see, type I unacceptable morphisms do always exist. 
We also discuss the more restrictive case of {\it discrete} and {\it stable} morphisms, 
which are quite meaningful from the point of view of the local Langlands correspondence.
We show that there is no type I discrete unacceptable morphisms, nor type II ones for $p>2$.
Moreover, although type III discrete unacceptable morphisms turn out to always exist, 
we show that there is no stable unacceptable discrete morphism for $p$ odd.  
All these results show that, for several natural families of Langlands parameters
for ${\rm PGSp}_6(F)$, the {\it weak} equivalence class appearing 
in the local Langlands correspondence in \cite[App. C, Thm. 12.6]{gansavin} 
coincides with the familiar, and stronger, one 
(given by conjugacy by the dual group).
\ps

In the final Sections~\ref{sec:gspinn} and~\ref{sec:noncompact},
we explain how our results can also be applied to study the acceptability of ${\rm GSpin}(n)$-valued 
morphisms, and as an example, we give an application to the acceptability of certain ${\rm GSpin}_7$-valued $\ell$-adic Galois representations which improves some result in \cite{kretshin}. 
\ps

\begin{remark} 
{\rm (A general remark on topology)}
\label{rem:topo}
In the study of unacceptable morphisms, we would not lose much 
in restricting to injective and continuous morphisms $\Gamma \rightarrow G$ from 
 compact groups $\Gamma$. Indeed, 
assume $r_1,r_2 : \Gamma \rightarrow G$ are two element conjugate morphisms, 
with $G$ a compact group and $\Gamma$ an arbitrary group.
Consider the morphism $r_1 \times r_2 : \Gamma \rightarrow G \times G$.
Up to replacing $\Gamma$ by its image under  
$r_1 \times r_2$,
we may assume $\Gamma \subset G \times G$ 
and that $r_1$ and $r_2$ are the two natural projections.
Define $\Gamma'$ as the (compact) closure of $\Gamma$ in $G \times G$, and
$r_1', r_2' : \Gamma' \rightarrow G$ as the two (continuous) projections. 
As $\{ (g,hgh^{-1}) \, | g,h \in G\}$ is closed in $G \times G$,
the morphisms $r'_1$ and $r'_2$ are element conjugate; 
they are conjugate in $G$ if, and only if, $r_1$ and $r_2$ are. 
In particular, for $(\gamma_1,\gamma_2) \in \Gamma'$ we have 
$\gamma_1=1$ if, and only if, $\gamma_2=1$, as $\gamma_1$ and $\gamma_2$
are conjugate, so $r_1$ and $r_2$ are both injective.

\end{remark}

\tableofcontents

\section{General notations on Spin groups}
\label{sec:defspin}

Let $n\geq 1$ be an integer and $E$ the standard Euclidean space $\R^n$
with inner product denoted by $x.y$ for $x,y \in E$. 
We denote by ${\rm O}(n)={\rm O}(E)$ the orthogonal group of $E$, 
by ${\rm SO}(n)= {\rm SO}(E)$ its special orthogonal group, 
and by ${\rm Cl}(E)$ the Clifford algebra of $E$. 
We have a natural inclusion $E \subset {\rm Cl}(E)$.
For $e \in E$ with $e.e=1$, we have $e^2 = 1$ in ${\rm Cl}(E)$ and
the conjugation by $e$ in ${\rm Cl}(E)$ preserves 
the subspace $E$ and induces the opposite of the Euclidean reflection of 
$E$ about $e$.
Our convention is that ${\rm Pin}(n)={\rm Pin}(E)$ is the subgroup of ${\rm Cl}(E)^{\times}$ 
generated by the elements $e \in E$ with $e.e=1$. 
The $\Z/2\Z$-grading of ${\rm Cl}(E)$ defines a group morphism 
${\rm deg} : {\rm Pin}(n) \rightarrow \Z/2\Z$ sending any such $e$ to $1 \in \Z/2\Z$, 
and whose kernel is by definition ${\rm Spin}(n)={\rm Spin}(E)$. 
We denote by $$\pi : {\rm Pin}(n) \rightarrow {\rm O}(n)$$ the group morphism 
defined for all $\gamma \in {\rm Pin}(n)$ and all $v \in E$ 
by the equality $$\pi(\gamma)(v) \,=\, (-1)^{\rm deg(\gamma)} \,\gamma v \gamma^{-1}$$
in ${\rm Cl}(E)$.  
The morphism $\pi$ is surjective as Euclidean reflections generate ${\rm O}(E)$. 
Its kernel is a central subgroup of order $2$, 
generated by an element of ${\rm Spin}(n)$ denoted $-1$. 
We have $\det \circ \pi = (-1)^{\rm deg}$ on ${\rm Pin}(n)$. 
Also, ${\rm Pin}(n)$ is a compact subgroup of the invertible elements 
of the finite dimensional $\R$-algebra ${\rm Cl}(E)$ and for this topology the morphism $\pi$ is continuous.
\ps

Assume now $n$ is odd and write $n=a+b$ with $a$ odd and $b \neq 0$ even. 
The $a$-dimensional subspaces of $E=\R^n$ form a single orbit under ${\rm SO}(E)$.
Fix such a subspace $A \subset E$ and set $B = A^\perp$.
The stabilizer in ${\rm SO}(E)$ of $A$, 
hence of $B$, is 
$${\rm S}({\rm O}(A) \times {\rm O}(B)) := \{ (g,h) \in {\rm O}(A) \times {\rm O}(B)\, \, |\, \, \det g = \det h\}.$$ 
It is isomorphic to ${\rm SO}(a) \times {\rm O}(b)$. 
We denote by ${\rm Spin}(A,B)$ the inverse image of ${\rm S}({\rm O}(A) \times {\rm O}(B))$ in ${\rm Spin}(E)$. 
As $B$ is nonzero, we have a natural order $2$ character ${\rm S}({\rm O}(A) \times {\rm O}(B)) \rightarrow \{ \pm 1\}$
sending $(g,h)$ to $\det g = \det h$. 
Composing this character with $\pi : {\rm Spin}(A,B) \rightarrow {\rm S}({\rm O}(A) \times {\rm O}(B))$ defines a character
\begin{equation}\label{eq:defkapa} \kappa : {\rm Spin}(A,B) \rightarrow \{ \pm 1\}. \end{equation}
The kernel of $\kappa$, the inverse image of ${\rm SO}(A) \times {\rm SO}(B)$ in ${\rm Spin}(E)$, 
coincides with $${\rm Spin}(A) \cdot {\rm Spin}(B) \simeq {\rm Spin}(A) \times {\rm Spin}(B)/ \langle (-1,-1) \rangle.$$ 
Let $e$ and $f$ be elements of $A$ and $B$ respectively, with $e . e = f . f = 1$. 
Then $\pi(ef)$ acts by $-{\rm id}$ on $\R e \perp \R f$, and by $+{\rm id}$ on its orthogonal.
So $ef$ is an element of ${\rm Spin}(A,B)$ with $\kappa(ef)=-1$. 
It satisfies $(ef)^2 = efef=-eeff=-1$.  
As $a$ is odd, the center of ${\rm Spin}(A)\cdot {\rm Spin}(B)$ is the center of ${\rm Spin}(B)$,
which has four elements since $b$ is even.  This center contains $-1$, but also the element 
\begin{equation} \label{defzb} z_B := f_1 f_2 \cdots f_b \in {\rm Spin}(B) \end{equation}
 where $f_i$ is some orthonormal basis of $B$,
 since $z_B$ anti-commutes with any $f_i$. Note that 
 $\pi(z_B)$ is the element $({\rm id}_A, -{\rm id}_B)$ of ${\rm S}({\rm O}(A) \times {\rm O}(B))$.
 It follows that $\pm z_B$ does not depend on the choice of the orthonormal basis $f_i$ of $B$.
 We also have $z_B^2 = (-1)^{b/2}$ and $z_B\, f\, = \,-\, f \,z_B$ for all $f \in B$.\ps
 
In the standard case $A=\R^a$, $B=\R^b$, $E= A \perp B = \R^{a+b}$,
 we simply write ${\rm Spin}(a,b)$ for ${\rm Spin}(A,B)$. 
 Also, by a {\it ${\rm Spin}(a,b)$-subgroup} of ${\rm Spin}(n)$ we mean a subgroup of the form ${\rm Spin}(A,B)$
 with $\dim A=a$ and $\dim B=b$; they form a single conjugacy class under ${\rm Spin}(n)$.
 
\begin{lemma} 
\label{lem:critacc} 
Assume we have $E= A \perp B$ with $a = \dim A$, $b= \dim B$, and $b>0$ even.
Define $\kappa : {\rm Spin}(A,B) \rightarrow \{ \pm 1\}$ as in Formula~\eqref{eq:defkapa}, and $\pm z_B \in {\rm Spin}(A,B)$
by Formula~\eqref{defzb}.  
Then for all $g \in {\rm Spin}(A,B)$ we have
$z_B \,g\, z_B^{-1} \,=\, \kappa(g) \,g$.
\end{lemma}

\begin{pf} For $\kappa(g)=1$, 
this follows as $z_B$ is in the center of $\ker \kappa = {\rm Spin}(A) \cdot {\rm Spin}(B)$. 
For $g=ef$ with $e \in A, f \in B$ and $e.e=f.f=1$, we have $z_B\,g\,z_B^{-1}\,=\,-g$ 
since $z_B$ commutes with $e$ (as $b$ is even) 
and anti-commutes with $f$ (as we have seen). 
We conclude as $ef$ and $\ker \kappa$ generate ${\rm Spin}(A,B)$. 
\end{pf}

\begin{remark}
\label{rem:centSpinAB}
Lemma~\ref{lem:critacc} implies that the center of ${\rm Spin}(A,B)$ is $\{ \pm 1\}$.
\end{remark}

\section{The conditions {\rm (U1)} and {\rm (U2)} for ${\rm Spin}(n)$ with general odd $n$}
\label{sect:generaln}

In all this section, $n \geq 1$ is an {\it odd} integer and $\Gamma$ is a group.

\begin{prop}
\label{defeta} \label{prop:defeta} 
Assume $r,r' : \Gamma \rightarrow {\rm Spin}(n)$ 
are element conjugate {\rm (}{\it i.e.} satisfy {\rm (U1)}{\rm )}. 
There is $g \in {\rm Spin}(n)$ and a character $\eta :  \Gamma \rightarrow \{ \pm 1\}$  such that 
for all $\gamma \in \Gamma$ we have $g \,r'(\gamma)\, g^{-1} \,=\,\eta(\gamma)\, r(\gamma)$.
\end{prop}

\begin{pf} 
As $n$ is odd, we know that the group ${\rm SO}(n)$ is acceptable. 
Since $\pi : {\rm Spin}(n) \rightarrow {\rm SO}(n)$ is surjective with kernel $\{\pm 1\}$, 
it follows that there is $g \in {\rm Spin}(n)$ such that 
for all $\gamma \in \Gamma$ there is $\eta(\gamma) \in \{ \pm 1\}$ 
with $g\, r(\gamma) \,g^{-1}\, =\, \eta(\gamma)\, r(\gamma)$. 
This defines a map $\eta : \Gamma \rightarrow \{ \pm 1\}$, which is necessarily a 
group morphism. 
\end{pf}

\begin{definition}
\label{def:retaU12}
Consider group homomorphisms 
$r : \Gamma \rightarrow {\rm Spin}(n)$ and $\eta : \Gamma \rightarrow \{\pm 1\}$. We 
say that $(r,\eta)$ satisfies {\rm (U1)} {\rm (}resp. {\rm (U2)}{\rm )} 
if this property holds with $r':=\eta r$. We say that $(r,\eta)$ is unacceptable if it satisfies both {\rm (U1)} and {\rm (U2)}. 
\end{definition}

Proposition \ref{defeta} asserts that for all unacceptable $r$ there is 
$\eta$ such that $(r,\eta)$ is unacceptable. 
Of course, if $(r,\eta)$ is unacceptable then we have $\eta \neq 1$. 
Another simple property is the following.

\begin{prop} 
\label{prop:facteta}
Let $r$ and $\eta$ be as in Definition~\ref{def:retaU12} 
and assume $(r,\eta)$ satisfies {\rm (U1)}.
Then $\eta$ is trivial on the kernel of the morphism
$\pi \circ r : \Gamma \rightarrow {\rm SO}(n)$.
\end{prop}

\begin{pf} Assume $\gamma \in \Gamma$ satisfies $\pi(r(\gamma))=1$, 
{\it i.e.} $r(\gamma) = \pm 1$. 
By (U1), $r(\gamma)$ is conjugate to $\eta(\gamma)r(\gamma)$ in ${\rm Spin}(n)$.
As $1$ and $-1$ are not conjugate, this forces $\eta(\gamma)=1$.
\end{pf}

\begin{definition}
\label{def:Xr}
Let $r : \Gamma \rightarrow {\rm Spin}(n)$ be a group morphism. 
We denote by ${\rm X}(r)$ the set of group morphisms $\chi : \Gamma \rightarrow \{ \pm 1\}$ such that 
there is $g \in {\rm Spin}(n)$ satisfying 
$$\forall \gamma \in \Gamma, \,\, \, g \,r(\gamma) \,g^{-1} = \chi(\gamma) r(\gamma).$$
This set ${\rm X}(r)$ is a subgroup of ${\rm Hom}(\Gamma,\{ \pm 1\})$.
\end{definition}

By definition, if $(r,\eta)$ satisfies (U1), then $(r,\eta)$ is unacceptable 
if, and only if, $\eta \notin {\rm X}(r)$. 

\begin{remark}
\label{rem:etar}
Let us denote by 
${\rm E}(r)$ the set of morphisms $\eta : \Gamma \rightarrow \{ \pm 1\}$ 
such that $(r,\eta)$ is unacceptable. By Proposition~\ref{prop:defeta}, $r$ is unacceptable if, and only if, 
${\rm E}(r)$ is nonempty. By Definition~\ref{def:Xr}, if we have $\eta \in {\rm E}(r)$
and $\chi \in {\rm X}(r)$, then $\eta \chi$ also belongs to ${\rm E}(r)$, so that ${\rm X}(r)$ acts freely by multiplication on ${\rm E}(r)$. 
\end{remark}

Our first aim now is to give an alternative description of ${\rm X}(r)$.
If $r : \Gamma \rightarrow {\rm Spin}(n)$ is a given morphism, 
it has a natural $n$-dimensional real representation 
$\pi \circ r : \Gamma \rightarrow {\rm SO}(n)$ on the Euclidean space $E=\R^n$.
We shall simply denote by $E$ this $\R[\Gamma]$-module.
This is a semi-simple $\R[\Gamma]$-module : 
for each $\R[\Gamma]$-submodule $V \subset E$, we have 
the $\Gamma$-stable decomposition $E = V \oplus V^\perp$.
For any such $V$, we also denote by $\det_V : \Gamma \rightarrow \{ \pm 1\}$ its determinant character.
We have $\det_E =1$, hence the equality $\det_V = \det_{V^\perp}$. \ps

\begin{prop} 
\label{prop:critU2}
Let $r : \Gamma \rightarrow {\rm Spin}(n)$ be a group morphism. 
The subgroup ${\rm X}(r)$ of ${\rm Hom}(\Gamma,\{ \pm 1\})$ is the subset of
characters of the form $\det_V$ where $V$ runs among the $\R[\Gamma]$-submodules of $E$.
\end{prop}

\begin{pf}
Assume we have an $\R[\Gamma]$-stable decomposition $E = A \perp B$.
Up to exchanging $A$ and $B$, 
we may assume $b:=\dim B$ is even and $>0$, and we set $a= n-b=\dim A$.
We have $r(\Gamma) \subset {\rm Spin}(A,B)$. 
The restriction to $\Gamma$ of the character 
$\kappa$ of Formula \eqref{eq:defkapa} is $\det_A=\det_B$ by construction.
Lemma~\ref{lem:critacc} thus shows that $\det_A$ and $\det_B$ are in ${\rm X}(r)$. \ps
In order to prove the proposition, it is enough to show that 
${\rm X}(r)$ is generated by the elements of the form $\det_V$ with 
$V \subset E$ an $\R[\Gamma]$-submodule. Indeed, this assertion implies
first that ${\rm X}(r)$ is generated by those $\det_V$ with $V$ irreducible,
and then using that these characters have order $\leq 2$, we deduce that
any element of ${\rm X}(r)$ has the form $\det_V$ with 
$V \subset E$ an $\R[\Gamma]$-submodule (non necessarily irreducible).\ps
\ps

Define ${\rm C}(r)$ as the centralizer of $\pi(r(\Gamma))$ in ${\rm SO}(n)$, 
and set ${\rm D}(r) = \pi^{-1}({\rm C}(r))$. In other words, we have
$${\rm D}(r) = \{ g \in {\rm Spin}(n)\,\,|\, \, \forall \gamma \in \Gamma, \, \, g \,r(\gamma)\, g^{-1} \,=\, \pm \,r(\gamma) \}.$$
This is a closed subgroup of ${\rm Spin}(n)$ containing $\{ \pm 1\}$, and that 
sits in the exact sequence 
$ 1 \rightarrow \{\pm 1\} \rightarrow {\rm D}(r) \overset{\pi}{\rightarrow} {\rm C}(r) \rightarrow 1$.
Fix $g \in {\rm D}(r)$. 
For each $\gamma \in \Gamma$, we have a unique sign 
${\rm e}_g(\gamma) \in \{\pm 1\}$ such that 
$$g \,r(\gamma)\, g^{-1} \,=\, {\rm e}_g(\gamma) r(\gamma).$$  
This formula shows that 
${\rm e}_g : \Gamma \rightarrow \{ \pm 1\}, \gamma \mapsto {\rm e}_g(\gamma)$, 
is a character, which clearly belongs to ${\rm X}(r)$, 
and also that ${\rm D}(r) \rightarrow {\rm X}(r), g \mapsto {\rm e}_g,$ 
is a group morphism. This latter morphism is surjective. 
Indeed, for $\chi \in {\rm X}(r)$
the element $g$ given by Definition~\ref{def:Xr} 
lies in ${\rm D}(r)$ and satisfies ${\rm e}_g=\chi$.
\ps

For a fixed $\gamma \in \Gamma$,
the morphism ${\rm D}(r) \rightarrow \{ \pm 1\}, g \mapsto {\rm e}_g(\gamma)$, 
is continuous, and trivial on the subgroup $\{ \pm 1\}$ of ${\rm D}(r)$, so it only depends
on the image of $g$ in ${\rm C}(r)/{\rm C}(r)^0$. 
Here, ${\rm C}(r)^0$ denotes the neutral connected component of ${\rm C}(r)$, 
so that we also have ${\rm C}(r)/{\rm C}(r)^0 \simeq \pi_0({\rm C}(r))$. 
As as a consequence,
we have proved that the (surjective) morphism 
${\rm D}(r) \rightarrow {\rm X}(r), g \mapsto {\rm e}_g$,
induces a (surjective) morphism 
$$\overline{{\rm e}} : {\rm C}(r)/{\rm C}(r)^0 \longrightarrow {\rm X}(r), $$
such that for all $g \in {\rm D}(r)$ we have $\overline{{\rm e}}( \pi (g)) = {\rm e}_g$. \ps

For each $\R[\Gamma]$-stable decomposition $E = A \oplus B$ 
as in the first paragraph of the proof, with $\dim B = b$ even and $>0$, 
we have $r(\Gamma) \subset {\rm Spin}(A,B)$ and Lemma~\ref{lem:critacc} 
precisely states that the element $z_B$ {\it loc. cit.} satisfies
\begin{equation} \label{eq:formZB} 
z_B \in {\rm D}(r)\, \, \,\,{\rm and}\, \,\,\, \overline{{\rm e}}(\pi(z_B)) \,=\, {\rm det}_B. 
\end{equation}
In order to conclude it remains to show that the $\pi(z_B)$, 
for $B$ an even dimensional $\R[\Gamma]$-submodule of $E$, do generate 
$\pi_0({\rm C}(r))$. Consider for this the isotypical decomposition 
$$E = \underset{i \in I}{\perp} E_i,$$
of the $\R[\Gamma]$-module $E$. 
The centralizer of $\pi \circ r$ in ${\rm O}(E)$ (rather than in ${\rm SO}(E)$)
is the direct product of the centralizer $C_i$ of $(\pi \circ r)_{|E_i}$ in 
${\rm O}(E_i)$. Write $E_i  \simeq U_i^{\oplus n_i}$ with $U_i$ an irreducible $\R[\Gamma]$-module 
and $n_i \geq 1$. 
According to Schur's lemma, there are three well-known possibilities: (see {\it e.g.} \cite[p. 96]{tdieck}) \ps

(a) ${\rm End}_{\R[\Gamma]}(U_i)\,=\,\R$. 
In this case, the $\C[\Gamma]$-module $U_i \otimes_\R \C$ is irreducible 
(we say that $U_i$ is {\it absolutely irreducible}), we have $C_i \simeq {\rm O}(n_i)$
and $(\det_{E_i})_{|C_i} = \det^{\dim U_i}$.\ps

(b) ${\rm End}_{\R[\Gamma]}(U_i) \simeq \C$ and $C_i \simeq {\rm U}(n_i)$ (unitary group). \ps

(c) ${\rm End}_{\R[\Gamma]}(U_i) \simeq \mathbb{H}$ and $C_i \simeq {\rm Sp}(n_i)$ (compact symplectic group). \ps

\noindent For $i \in I$ of type (b) or (c) the group $C_i$ is connected, so we have 
$C_i \subset {\rm SO}(E)$ and $C_i \subset {\rm C}(r)^0$.
It follows that we have $${\rm C}(r)/{\rm C}(r)^0 \simeq (\Z/2\Z)^m$$ 
where $m$ or $m+1$ is the number
of $i \in I$ of type $(a)$. More precisely, choose for each $i \in I$ of type $(a)$, 
an $\R[\Gamma]$-stable decomposition 
$E_i = F_i \perp G_i$ with $F_i$ irreducible (so $F_i \simeq U_i$),
and consider the element $\sigma_i \in {\rm O}(E)$ acting by $-{\rm id}$ on $F_i$ 
and by ${\rm id}$ on $G_i$ and each $E_j$ with $j \neq i$. 
This element is in $C_i \smallsetminus C_i^0$, and has determinant $(-1)^{\dim U_i}$.
It follows that ${\rm C}(r)/{\rm C}(r)^0$ is generated by 
the images of the $\sigma_i$ with $\dim U_i$ even, 
and by the images 
of the $\sigma_i\sigma_j$ with $i\neq j$ and $\dim U_i \equiv \dim U_j\equiv 1\bmod 2$. 
Each of these elements has the form $\pi(z_B)$ for $B=F_i$ or $B=F_i \perp F_j$.
\end{pf} 

\begin{coro} 
\label{cor:anyXdetV}
${\rm X}(r)$ is generated by those $\det_V$ with $V$ irreducible.
\end{coro}

\begin{pf} Immediate from Proposition~\ref{prop:critU2}.\end{pf}

\begin{coro} 
\label{cor:coroetastable} 
Assume we have morphisms $r : \Gamma \rightarrow {\rm Spin}(n)$ 
and $\eta : \Gamma \rightarrow \{\pm 1\}$. Then $(r,\eta)$ satisfies {\rm (U2)} if, and only if,
there is no $\R[\Gamma]$-stable subspace $V$ of $E$ with $\det_V = \eta$.
\end{coro}

\begin{pf} 
Clear from the definition of ${\rm X}(r)$ and Proposition~\ref{prop:critU2}. 
\end{pf}

We will say that a group embedding ${\rm Spin}(n) \rightarrow {\rm Spin}(m)$, with $m\geq n$, is
{\it standard}  if the restriction to ${\rm Spin}(n)$ of the standard representation of ${\rm Spin}(m)$
is isomorphic to $E \oplus 1^{m-n}$ (with $1$ the trivial representation). 

\begin{coro}
\label{cor:corU2stab}
Assume $r : \Gamma \rightarrow {\rm Spin}(n)$ is unacceptable (recall $n$ is odd).
Let $m\geq n$ be an integer and $\rho : {\rm Spin}(n) \rightarrow {\rm Spin}(m)$ a standard 
embedding. Then $\rho \circ r :  \Gamma \rightarrow {\rm Spin}(m)$ is unacceptable.
\end{coro}

\begin{pf} 
By Proposition~\ref{prop:defeta}, we may choose $\eta$ such that $(r,\eta)$ is unacceptable.
Note that have $\rho \circ (\eta r) \,=\, \eta\, (\rho \circ r)$.
Its thus clear that $(\rho \circ r,\eta)$ satisfies (U1). 
In order to check that it satisfies (U2), 
we may of course increase $m$, hence assume $m$ odd as well.  
We conclude by Corollary~\ref{cor:coroetastable} applied to $(r,\eta)$ and $(\rho \circ r,\eta)$.
\end{pf}

We now give a few equivalent conditions for (U1). 
Recall that we denote by $E$ the standard 
representation of ${\rm Spin}(n)$ (an $n$-dimensional real representation)
and by $S$ a Spin representation (a $2^{(n-1)/2}$-dimensional complex vector space).
If $r : \Gamma \rightarrow {\rm Spin}(n)$ is a given morphism, 
it allows to view $E$ as an $\R[\Gamma]$-module
and $S$ as a $\C[\Gamma]$-module.\ps

\begin{lemma} 
\label{lem:critamb} 
Assume $\gamma$ is in ${\rm Spin}(n)$ with $n$ odd. 
The following are equivalent:
\ps 
\begin{itemize}
\item[(i)] $\gamma$ is conjugate to $-\gamma$ in ${\rm Spin}(n)$, \ps\ps
\item[(ii)] $\gamma$ admits the eigenvalue $-1$ on $E$,\ps\ps
\item[(iii)] the trace of $\gamma$ is $0$ on $S$.\ps\ps
\end{itemize}
\end{lemma}

\begin{pf}
We first show the equivalence between (i) and (ii).
This could be proved directly, but we rather deduce it from Proposition~\ref{prop:critU2}.
Consider $\Gamma = \Z$ and $r : \Gamma \rightarrow {\rm Spin}(n)$ 
sending $m \in \Z$ to $\gamma^m$.
Define $\eta : \Z \rightarrow \{ \pm 1\}$ by $\eta(m) = (-1)^m$. 
Since $\Z$ is generated by $1$, assertion (i) holds if and only if we have $\eta \in {\rm X}(r)$.
By Proposition~\ref{prop:critU2} this is equivalent to ask that 
$E$ has a $\gamma$-stable subspace $V$ with $\det_V(\gamma)=-1$,
{\it i.e.} that $\gamma$ has the eigenvalue $-1$ in $E$.\ps

Let us now prove the equivalence between (i) and (iii). 
Note that the representation ring of the simply connected group ${\rm Spin}(n)$ 
is generated by its fundamental representations (over $\C$), 
which are $S$ and some others which all factor through ${\rm SO}(n)$ (since $n$ is odd).
It follows that two elements $\gamma_1$ and $\gamma_2$ 
of ${\rm Spin}(n)$ are conjugate  if, and only if, 
(a) their image are conjugate in ${\rm SO}(n)$, 
and (b) $\gamma_1$ and $\gamma_2$ have the same trace in $S$.
We conclude by applying this remark 
to $\gamma_1=\gamma$ and $\gamma_2=-\gamma$, 
since $-1$ acts by $-{\rm id}$ on $S$.
\end{pf}

\begin{remark} 
\label{rem:corprop}
For all $n$, an element $\gamma$ in ${\rm Spin}(n)$ 
is conjugate to $-\gamma$ if and only if $\gamma$ 
admits the eigenvalues $1$ and $-1$ on $E$.  
\end{remark}

We now give a few equivalent properties to (U1). 
For $r : \Gamma \rightarrow {\rm Spin}(n)$, 
we have already defined the $\R[\Gamma]$-module $E$ with $\dim_\R E=n$.
We also define the $\R[\Gamma]$-module 
\begin{equation}\label{def:LEsharp}
 \Lambda^{\sharp} E = \oplus_{i=1}^k \Lambda^i E, \, \, \, {\rm with}\, \,\, n=2k+1.
 \end{equation} 
The $\Gamma$-modules $E, S$ and $\Lambda^{\sharp} E$ are semi-simple, 
since they extend to the compact group ${\rm Spin}(n)$. 
Both $E$ and $\Lambda^{\sharp} E$ are defined over $\R$, 
and $S$ is defined over $\C$.

\begin{prop} 
\label{prop:propequat} 
Assume $r : \Gamma \rightarrow {\rm Spin}(n)$ is a group morphism (recall $n$ is odd)
and let $\eta : \Gamma \rightarrow \{ \pm 1\}$ be a morphism. 
The following are equivalent: \begin{itemize}\ps
\item[(i)] $(r,\eta)$ satisfies {\rm (U1)}, \ps\ps
\item[(ii)] for all $\gamma \in \Gamma$, $\gamma$ has the eigenvalue $\eta(\gamma)$ on $E$,\ps\ps
\item[(iii)] there is an isomorphism $\Lambda^{\sharp} E \simeq \Lambda^{\sharp} E \otimes \eta$ of $\R[\Gamma]$-modules,\ps\ps
\item[(iv)] there is an isomorphism $S \,\simeq \,\eta \otimes S$ of $\C[\Gamma]$-modules.\ps\ps
\end{itemize}
\end{prop}

\begin{pf} 
Since $1$ is an eigenvalue of each element of ${\rm SO}(n)$ for odd $n$,
the equivalence of (i) and (ii) follows from that of Lemma~\ref{lem:critamb}.
But (ii) is equivalent to $\det( \eta(\gamma) - r(\gamma)) =0$ 
for all $\gamma$ in $\Gamma$, or equivalently, 
\begin{equation} \label{cayleyham} 
\sum_{i=0}^{n}\,\,(-1)^i\,\, {\rm trace}( \gamma\, | \, \Lambda^i E)\,\, \eta(\gamma) ^{n-i} = 0
\end{equation} 
for all $\gamma$ in $\Gamma$. 
Set $X_+= \bigoplus_{i\,\, {\rm even}} \Lambda^i E$ 
and $X_{-}= \bigoplus_{i\,\,{\rm odd}} \Lambda^i E$. 
As an $\R[{\rm SO}(n)]$-module we have 
$E \simeq E^\star$ and $\det E = 1$, 
so $\Lambda^i E \simeq \Lambda^{n-i} E$ for all $i$. 
This shows $X_{+} \simeq X_{-} \simeq \Lambda^{\sharp} E$ as $n$ is odd. 
As $\Lambda^{\sharp} E$ is semi-simple, 
\eqref{cayleyham} is thus equivalent to 
$\Lambda^{\sharp} E \simeq \Lambda^{\sharp} E \otimes \eta$, 
and we have proved $(ii) \iff (iii)$. \par

We end by proving $(i) \iff (iv)$. Set $r' = \eta r$. 
The equivalence (i) $\iff$ (iii) of Lemma~\ref{lem:critamb} 
shows that $r$ and $r'$ are element conjugate if, and only if, 
$S$ and $\eta\otimes S$ have the same trace on $\Gamma$.
We conclude by semi-simplicity of the $\C[\Gamma]$-module $S$.
\end{pf}

\begin{remark}
Condition (iv) also means that we have 
$S \simeq {\rm Ind}_{\Gamma'}^{\Gamma} S'$ for some complex representation
$S'$ of the index $2$ subgroup $\Gamma' = \ker \eta$ of $\Gamma$.\ps
\end{remark}

Of course, Corollary~\ref{cor:coroetastable} and Proposition~\ref{prop:propequat} 
together imply Theorem~\ref{thm:gen} of the introduction.

\section{The case $n=7$}
\label{sect:n7}

In this section we focus on the case $n=7$, the first for which there turns out to exist unacceptable 
morphisms $\Gamma \rightarrow {\rm Spin}(n)$. 
Several arguments below follow Larsen's original arguments in \cite{larsen2}, 
and correct the erroneous Lemma 2.3 and Proposition 2.4 {\it loc. cit}.
An important role will be played by 
the conjugacy class of ${\rm Spin}(1,6)$-subgroups of ${\rm Spin}(7)$.\ps\ps

{\bf Notation:} {\it We fix a decomposition $\R^7 = L \perp F$ with $L$ a line, 
and we denote by $N$ the associated ${\rm Spin}(1,6)$-subgroup 
${\rm Spin}(L,F)$ of ${\rm Spin}(7)$.} \ps \ps

We have a natural surjective morphism 
$N \overset{\pi}{\rightarrow} {\rm S}({\rm O}(L) \times {\rm O}(F)) \simeq {\rm O}(F)$.
The inverse image of ${\rm SO}(F) \simeq {\rm SO}(6)$ in $N$ is isomorphic to ${\rm Spin}(6)$, 
has index $2$ in $N$, hence coincides with the connected component $N^0$ of the identity in $N$.
The choice of a unitary half-spin representation of $N^0 \simeq {\rm Spin}(6)$ 
identifies it with ${\rm SU}(4)$. 
We definitely fix such an identification and allow 
ourselves to write $N^0={\rm SU}(4)$ accordingly. 
The precise structure of $N$ is the following.

\begin{lemma} 
\label{lem:strucN}
The group $N$ is the semi-direct product of $\Z/2\Z$ by its subgroup ${\rm SU}(4)$
with respect to an order $2$ symmetric outer automorphism. 
\end{lemma}

The meaning of the statement is the following.
Recall that for $n\geq 2$ {\it even}, any order $2$ outer automorphism 
$\theta$ of ${\rm SU}(n)$ has the form
$g \mapsto p \overline{g} p^{-1}$ with $p \in {\rm SU}(n)$ satisfying ${}^{\rm t}\!p = \pm p$; 
we say that $\theta$ is {\it symmetric} if $p$ is, 
and {\it antisymmetric} otherwise. 
The subgroup of fixed points in ${\rm SU}(n)$ of such a $\theta$ 
is isomorphic to ${\rm SO}(n)$ in the symmetric case, and to ${\rm Sp}(n/2)$ in the antisymmetric case. 
Also, in the natural semi-direct product 
${\rm SU}(n) \rtimes \langle \theta \rangle$,
the elements of the form $\vartheta=h \theta$ with $h \in {\rm SU}(n)$ 
and $\vartheta^2=1$ (resp. $\vartheta^2=-1$)
induce by conjugation all the order $2$ outer automorphisms of ${\rm SU}(n)$ with same 
(resp. opposite) symmetry type as $\theta$.  
In particular, the lemma above asserts that there is $\vartheta \in N \smallsetminus N^0$ 
with $\vartheta^2=1$ and $\vartheta g \vartheta^{-1} = \overline{g}$ for all $g \in {\rm SU}(4)$.

\begin{pf} 
Write $L = \R e$ and choose an orthonormal basis $f_1,\dots,f_6$ of $F$.
Consider the element $\tau = e f_1 f_2 f_3 \in {\rm Spin}(7)$, which satisfies $\tau^2=1$.
We have $\tau \in N \smallsetminus N^0$, as $\pi(\tau)$ acts by $-1$ on $L$, and by a symmetry on $F$
with fixed subspace of dimension $3$. 
This shows $N = N^0 \rtimes \langle \tau \rangle$. 
The conjugation by $\tau$ defines an automorphism of $N^0$;
it is not inner
since we have $\tau z_F \tau^{-1} = -z_F$ by Lemma~\ref{lem:critacc},
whereas $z_F$ is central (of order $4$) in $N^0$.
The morphism $\pi$ induces an exact sequence on fixed points 
$1 \rightarrow \{ \pm 1\} \rightarrow {\rm SU}(4)^{\tau=1} \rightarrow {\rm O}(6)^{\pi(\tau)=1}$.
But the description above of $\pi(\tau)$ shows 
that the neutral component of ${\rm O}(6)^{\pi(\tau)=1}$
is ${\rm SO}(3) \times {\rm SO}(3)$. 
The only possibility is thus ${\rm SU}(4)^{\tau=1} \simeq {\rm SO}(4)$, 
and we are done.
\end{pf}

\begin{remark} 
\label{rem:pin6ornotpin6}
A similar argument shows that ${\rm Pin}(6)$ is also 
a semi-direct product of $\Z/2\Z$ by ${\rm SU}(4)$, 
but with respect to an antisymmetric outer automorphism. 
In particular, it is not isomorphic to ${\rm Spin}(1,6)$.\footnote{As an example, 
the compact form of the Langlands dual group of 
${\rm PU}(4)$ is isomorphic to ${\rm Pin}(6)$, and not to ${\rm Spin}(1,6)$.}
As the embeddings ${\rm Spin}(6) \rightarrow {\rm Spin}(7)$ are unique 
up to ${\rm Spin}(7)$-conjugacy,
and with normalizers the ${\rm Spin}(1,6)$-subgroups, 
it follows that the compact group ${\rm Pin}(6)$ does not embed into ${\rm Spin}(7)$.
\end{remark}
\ps
We will now give a few properties of the $N$-valued morphisms.  
\ps
\begin{remark} 
\label{rem:conjinsN}
Let $G$ be a semi-direct product of $\Z/2\Z$ by ${\rm SU}(m)$ defined by an outer automorphism.
Two morphisms $r,r' : \Gamma \rightarrow {\rm SU}(m)$ are conjugate in $G$ if, and only if, 
we have $r \simeq r'$ or $r^\ast \simeq r'$ as 
$m$-dimensional representations of $\Gamma$. 
\end{remark}

Recall that $\kappa : N \rightarrow \{ \pm 1\}$ 
denotes the character of $N$ acting on the line $L \subset E$, 
or equivalently, the determinant of the action of $N$ on $F \subset E$  (see Formula \eqref{eq:defkapa}).
We have\, $\ker \kappa = N^0$.
The embedding $N \subset {\rm Spin}(7)$ has the following nice property.

\begin{prop} \label{prop:Nfus}
Let $r,r' : \Gamma \rightarrow N$ be two group morphisms.
The following are equivalent:\begin{itemize}\ps
\item[(i)] $r$ and $r'$ are conjugate under $N$,\ps
\item[(ii)] $r$ and $r'$ are conjugate in ${\rm Spin}(7) \supset N$, 
and we have $\kappa \circ r =\kappa \circ r'$.\ps
\end{itemize}
\end{prop}

In particular, two elements of $N$ are conjugate if, and only if, they have the same image in 
$N/N^0$ and are conjugate in ${\rm Spin}(7)$ (case $\Gamma=\Z$).

\begin{pf} 
We clearly have (i) $\implies$ (ii).
Assume (ii) holds.  
In particular, the two $7$-dimensional representations $\pi \circ r$ 
and $\pi \circ r'$ of $\Gamma$ on $E$ are isomorphic.
They both have the same character on the line $L$, 
as this character is $\kappa \circ r = \kappa \circ {r'}$, 
hence their restriction to $F=L^\perp$ are isomorphic as well, 
hence ${\rm O}(F)$-conjugate by Lemma~\ref{lem:straccclassical}.
Up to replacing $r'$ by an $N$-conjugate, 
we may thus assume that we have $\pi \circ r = \pi \circ r'$, 
{\it i.e. } $r' = \eta r$ for some morphism $\eta : \Gamma \rightarrow \{\pm 1\}$. 
By definition and (ii), we have then $\eta \in {\rm X}(r)$. 
Consider an irreducible decomposition $F= \oplus_i  F_i$ of  
the $\R[r(\Gamma)]$-submodule of $F$. 
Note that the elements $\pm z_F$ are in $N$, 
as well as the $\pm z_{F_i}$ for $\dim F_i$ even, 
and the $z_L z_{F_i}$ for $F_i$ odd; 
with respective image in ${\rm Hom}(\Gamma,\{\pm 1\})$ 
the characters $\det_L = \kappa \circ r$, 
$\det_{F_i}$ and $\det_L \det_{F_i}$, by Formula~\eqref{eq:formZB}. 
Since ${\rm X}(r)$ is generated by $\det_L$ 
and the $\det_{F_i}$ by Corollary~\ref{cor:anyXdetV}, 
it follows that there is $n \in N$ such that $n r n^{-1} = \eta  r$, 
which proves (i).
\end{pf}

We have an exact sequence $1 \rightarrow \{ \pm 1\} \rightarrow N \rightarrow {\rm O}(F) \rightarrow 1$.
As $\{\pm 1\}$ is central in $N$, and as ${\rm O}(F)$ is acceptable, we have the:

\begin{propdefi}
\label{def:defetaN} The statement of Proposition~\ref{defeta} also holds with ${\rm Spin}(n)$
replaced by $N$, with the same proof. Similarly, Definition~\ref{def:retaU12} also makes sense
for pairs of morphisms $r: \Gamma \rightarrow N$ and $\eta : \Gamma \rightarrow \{ \pm 1\}$.
\end{propdefi}

This definition does not conflict with Definition~\ref{def:retaU12}, more precisely:

\begin{coro} 
\label{coro:embNSpin7}
Let $r : \Gamma \rightarrow N$ and $\eta : \Gamma \rightarrow \{\pm 1\}$ be two morphisms,
and $\widetilde{r} : \Gamma \rightarrow {\rm Spin}(7)$ the composition of $r$ with the inclusion $N \subset {\rm Spin}(7)$. Then $(r,\eta)$ satisfies {\rm (U1)} {\rm (}resp. {\rm (U2)}{\rm )} if and only if $(\widetilde{r},\eta)$ has this property. In particular, $(r,\eta)$ is unacceptable if, and only if, $(\widetilde{r},\eta)$ is unacceptable.
\end{coro}

\begin{pf} 
If $r$ and $r'=\eta r$ are ${\rm Spin}(7)$-conjugate, they are $N$-conjugate 
by the proposition, as we have $\kappa \circ r' = \kappa \circ r$ as $-1 \in N^0$. 
In the special case $\Gamma=\Z$, this also shows that if $\widetilde{r}$ and $\eta \widetilde{r}$
are element conjugate in ${\rm Spin}(7)$, then $r$ and $\eta r$ are element conjugate in $N$.
The reverse implications are trivial.
\end{pf}

\begin{prop}
\label{prop:pin16unacc}
The group ${\rm Spin}(1,6)$ {\rm (}hence $N${\rm )} is not acceptable. 
\end{prop}


We give below a first example of an unacceptable $N$-valued morphism.
We denote by $\mu_d \subset \C^\times$ the subgroup of $d$-th roots of unity.

\ps
\medskip
\begin{center} {\sc Example 1}\end{center}
\ps

Set $\Gamma=\mu_4 \times \mu_4$ and denote by $a$ and $b$ the two order $4$ characters of $\Gamma$ 
defined by the first and second projections respectively. 
Consider a morphism $r : \Gamma \rightarrow {\rm SU}(4)$ which, 
viewed as a complex $4$-dimensional representation of $\Gamma$, satisfies
$$r  \simeq a \oplus ab^2 \oplus b \oplus ba^2.$$
Such a morphism $r$ is unique up to ${\rm SU}(4)$-conjugation by Lemma~\ref{lem:straccclassical}.
Set $\eta = b^2$. We have in particular 
$\eta r  \simeq  r \otimes \eta \simeq a \oplus ab^2 \oplus b^{-1} \oplus b^{-1}a^2$. 
As neither $b^{-1}$ nor $a^{-1}$ appears in $r$, observe that neither $\eta r$, 
nor its dual, is isomorphic to $r$. 
This shows that $r$ and $\eta r$ are not $N$-conjugate, by Lemma~\ref{lem:strucN} and Remark~\ref{rem:conjinsN}.
Nevertheless, $r$ and $\eta r$ are element conjugate:  for all $g \in \Gamma$, we have $\eta(g) r(g)$ conjugate in ${\rm SU}(4)$ to either $r(g)$ or $r(g)^{-1}$. 
Indeed, this is obvious for $g$ in ${\rm ker} \,\eta$. 
Moreover, over any of the two subgroups ${\rm ker}\,\,a^2$ and ${\rm ker}\, \,(ab)^2$ we have 
$a \oplus ab^{2} \simeq a^{-1}\oplus a^{-1}b^2$, so $r(g)^{-1}$ is conjugate to $\eta(g) r(g)$. 
We conclude as $\Gamma$ is the union of ${\rm ker}\, a^2$, ${\rm ker}\, b^2$ and ${\rm ker}\, (ab)^2$. 
$\square$
\ps

\begin{coro} 
\label{cor:spin7unacc}
For $n\geq 7$, the group ${\rm Spin}(n)$ is not acceptable.
\end{coro}

\begin{pf} 
For $n=7$, this follows from Example 1 and Corollary~\ref{coro:embNSpin7}.
The general case follows then from the case $n=7$ and Corollary~\ref{cor:corU2stab}.
\end{pf}

\begin{remark} 
\label{rem:repNE}
The representation of $N^0={\rm SU}(4)$ on $E \otimes \C$ via $\pi$ is 
$1 \oplus \Lambda^2 V$ where $V=\C^4$ is its tautological representation. As a consequence,
if $r : \Gamma \rightarrow {\rm SU}(4) \rightarrow {\rm Spin}(7)$ is as in Example 1, 
then we have a $\C[\Gamma]$-module isomorphism 
$$E \otimes \C \,\simeq\,  1 \oplus a^2b^2 \oplus a^2 b^2 \oplus \left( ab \oplus a^{-1}b^{-1} \right) \oplus \left( a^{-1} b \oplus ab^{-1} \right),$$ 
as well as ${\rm X}(r)=\{1,a^2b^2\}$ and $\eta = b^2$.
\end{remark}

Before giving more examples, our goal up to the end of this section will be to prove 
Theorem~\ref{thm:exceptionPin6} of the introduction, which is a key result of this paper.
Our proof below is inspired by Larsen's analysis in~\cite{larsen2}.
We start with some general facts about ${\rm Spin}(7)$.
Recall that the spin representation of ${\rm Spin}(7)$ is well-known to be irreducible, real, 
and the unique $8$-dimensional faithful representation of ${\rm Spin}(7)$.
By Lemma~\ref{lem:straccclassical}, there is thus a unique ${\rm O}(8)$-conjugacy class 
of embeddings ${\rm Spin}(7) \rightarrow {\rm SO}(8)$.
As every automorphism of ${\rm Spin}(7)$ is inner we obtain:

\begin{lemma} 
\label{lem:subso8}
There are exactly two conjugacy classes of compact subgroups of ${\rm SO}(8)$ 
which are isomorphic to ${\rm Spin}(7)$, and these two classes are 
permuted transitively under ${\rm O}(8)$-conjugation.
\end{lemma}

The center of each ${\rm Spin}(7)$ in ${\rm SO}(8)$ is the center $\{ \pm 1\}$ of ${\rm SO}(8)$. 
We denote by $Z$ the center of ${\rm Spin}(8)$; we have $Z \simeq \Z/2\Z \times \Z/2\Z$.
Recall that any {\it triality} automorphism of ${\rm Spin}(8)$ permutes transitively the three 
order $2$ subgroups of $Z$.

\begin{lemma} 
\label{lem:subspin8}
There are three conjugacy classes of compact subgroups of ${\rm Spin}(8)$ 
which are isomorphic to ${\rm Spin}(7)$. These three classes are 
distinguished by their intersection with $Z$, 
which can be any of the three order $2$ subgroups of $Z$.
\end{lemma}

\begin{pf} 
We denote by $\mathcal{A}$ the set of subgroups of ${\rm Spin}(8)$ which are isomorphic to ${\rm Spin}(7)$, and by $\mathcal{B}$ the set of subgroups of ${\rm SO}(8)$ which are isomorphic to either ${\rm Spin}(7)$ or ${\rm SO}(7)$.
Fix also a surjective morphism $\rho : {\rm Spin}(8) \rightarrow {\rm SO}(8)$. 
For $T$ in $\mathcal{A}$ we have $\rho(T) \simeq {\rm Spin}(7)$ if\, $(\ker \rho) \cap T =\{1\}$\, 
and $\rho(T) \simeq {\rm SO}(7)$ otherwise, hence $\rho$ induces a natural map $\rho_\ast : \mathcal{A} \rightarrow \mathcal{B}, T \mapsto \rho(T)$.  We claim that $\rho_\ast$ is bijective.\ps 
Indeed, the injectivity is clear over the $T \in \mathcal{A}$ mapped to ${\rm SO}(7)$, since then we have $T=\rho^{-1}(\rho(T))$,  and follows from ${\rm Hom}({\rm Spin}(7),\mu_2)=\{1\}$ for the others. 
For the surjectivity, note that as ${\rm Spin}(7)$ is simply connected, any $S \in \mathcal{B}$ 
which is isomorphic to ${\rm Spin}(7)$ has the form $\rho(T)$ 
for some (unique) $T \in \mathcal{A}$. To conclude, we recall the easy 
fact that the subgroups of ${\rm SO}(8)$ isomorphic to ${\rm SO}(7)$ 
are all conjugate, since they are the stabilizers of the norm $1$ elements of $\R^8$, 
and that their inverse image in ${\rm Spin}(8)$ are isomorphic to ${\rm Spin}(7)$. 
So $\rho_\ast$ is bijective. \ps

By Lemma~\ref{lem:subso8} and the previous sentence, 
there are three ${\rm SO}(8)$-conjugacy classes in $\mathcal{B}$. 
As $\rho$ is surjective, and as $\rho_\ast$ is bijective and commutes with conjugacy, 
there are also three ${\rm Spin}(8)$-conjugacy classes in $\mathcal{A}$, 
and a single one with center ${\rm ker}\, \rho$. 
We conclude as any order $2$ subgroup of $Z$ is the kernel of a suitable $\rho$, by {\it triality}.
\end{pf}
\ps\ps

As is well-known, the exceptional compact group ${\rm G}_2$ has a unique non trivial
irreducible represesentation of dimension $\leq 8$, and it is faithful of dimension $7$.
By Lemma~\ref{lem:straccclassical}, it follows that both ${\rm SO}(7)$ and ${\rm SO}(8)$ have a unique
conjugacy class of compact subgroups isomorphic to ${\rm G}_2$. Using that 
${\rm G}_2$ is simply connected and ${\rm Hom}({\rm G}_2,\mu_2)=\{1\}$, it follows 
that ${\rm Spin}(7)$ and ${\rm Spin}(8)$ also have a unique conjugacy class of compact subgroups isomorphic to ${\rm G}_2$, and that we shall call the {\it ${\rm G}_2$-subgroups}. 
The centralizer in ${\rm Spin}(7)$ of a ${\rm G}_2$-subgroup $H$ is thus the center $\{ \pm 1\}$ of ${\rm Spin}(7)$, and since we have ${\rm Out}({\rm G}_2)=1$, the normalizer of $H$ in ${\rm Spin}(7)$ is $\{ \pm 1\} \times H$. 
The following proposition corrects \cite[Prop. 2.4]{larsen2}.
\begin{prop} 
\label{prop:subspin7}
Let $S_1$ and $S_2$ be two different subgroups of ${\rm SO}(8)$ 
both isomorphic to ${\rm Spin}(7)$. 
Then $S_1 \cap S_2$ is one of the following subgroups of
\footnote{The choice of the isomorphism  $S_1 \simeq {\rm Spin}(7)$ does not matter as ${\rm Out}({\rm Spin}(7))=1$. } $S_1 \simeq {\rm Spin}(7)$:
\begin{itemize}
\item[(i)] the normalizer of a ${\rm G}_2$-subgroup,\ps
\item[(ii)] a ${\rm Spin}(1,6)$-subgroup,\ps
\item[(iii)] the identity component of a ${\rm Spin}(1,6)$-subgroup.\ps
\end{itemize}
We are in the first case if, and only if, $S_1$ and $S_2$ are not conjugate under ${\rm SO}(8)$.
\end{prop}

\begin{pf}  
Fix a surjective morphism $\pi : {\rm Spin}(8) \rightarrow {\rm SO}(8)$ 
and set $\widetilde{S}_i = \pi^{-1} S_i$ for $i=1,2$. 
By the bijectivity of $\rho_\ast$ in the previous proof, there are unique subgroups 
$T_i \subset {\rm Spin}(8)$ isomorphic to ${\rm Spin}(7)$ 
and with $\pi(T_i)=S_i$. Moreover, we have 
$T_i\, \cap \,(\ker \,\pi) \,=\,1$, $\widetilde{S}_i \,=\, Z T_i$, 
and the groups $S_i$ are conjugate in ${\rm SO}(8)$ if, 
and only if, the centers of $T_1$ and $T_2$ coincide. 
By triality, we may choose a morphism $\rho : {\rm Spin}(8) \rightarrow {\rm SO}(8)$ 
whose kernel is the center of $T_1$ (so $\rho \neq \pi$). 
The subgroup $\rho(\widetilde{S}_1)=\pm \rho(T_1)$ of 
${\rm SO}(8)$ is thus isomorphic to ${\rm O}(7) = \Z/2 \times {\rm SO}(7)$, 
namely it is the stabilizer of a line $L_1$ in the Euclidean $\R^8$. 
There are two cases: \ps

(a) {\it $S_1$ and $S_2$ are not conjugate under ${\rm SO}(8)$}. 
We have then $\rho(T_2) \simeq {\rm Spin}(7)$ and 
the subgroup $\rho(\widetilde{S}_2)=\pm \rho(T_2)=\rho(T_2)$ of ${\rm SO}(8)$ 
is thus also isomorphic to ${\rm Spin}(7)$. But the subgroups 
of ${\rm Spin}(7)$ isomorphic to ${\rm G}_2$ (resp. $\Z/2\Z \times G_2$) 
are exactly the stabilizers of nonzero vectors (resp. lines) 
in the spin representation of ${\rm Spin}(7)$ by \cite[Thm. 5.5]{adams}.  
Applying this to $\rho(\widetilde{S}_2)$ we obtain
$$\rho(\widetilde{S}_1\cap \widetilde{S}_2) = \rho(\widetilde{S}_1) \cap \rho(\widetilde{S}_2) \simeq \Z/2\Z \times {\rm G}_2.$$
As $G_2$ is simply connected and satisfies ${\rm Hom}(G_2,\mu_2)=1$, 
we obtain $\pi^{-1} (S_1 \cap S_2)  =\widetilde{S}_1\cap \widetilde{S}_2= Z \times H$ 
with $H \simeq {\rm G}_2$, hence $S_1 \cap S_2 \simeq \Z/2\Z \times {\rm G}_2$. This is case (i). \ps

(b) {\it $S_1$ and $S_2$ are conjugate under ${\rm SO}(8)$}. 
In this case $T_1$ and $T_2$ have the same center, namely $\ker \, \rho$.
The subgroup $\rho(\widetilde{S}_2)$ is then isomorphic to ${\rm O}(7)$ as well, 
{\it i.e.} it is the stabilizer of a line $L_2$ in $\R^8$. 
We have $L_1 \neq L_2$ as $S_1 \neq S_2$. 
Let $Q$ be the Euclidean plane $L_1 \oplus L_2$ in $\R^8$, 
$U \subset {\rm O}(Q)$ the subgroup fixing pointwise $L_1$ and preserving $L_2$,
and $I \subset {\rm SO}(8)$ the subgroup of elements 
$(g,h)$ of $U \times {\rm O}(Q^\perp)$ with $\det g = \det h$. 
We have $\rho(\widetilde{S}_1) \cap \rho(\widetilde{S}_2) = \{ \pm 1\} \times I$ in ${\rm SO}(8)$. 
As $U$ preserves the orthogonal $L_3$ of $L_1$ in $Q$ (a ``third'' line in $Q$), 
there are two sub-cases:\ps
(b1)  $L_3=L_2$,  $U = {\rm SO}(L_1) \times {\rm O}(L_2) \simeq \Z/2\Z$ and $I \simeq {\rm O}(Q^\perp) \simeq {\rm O}(6)$, or \ps
(b2) $L_3 \neq L_2$,  $U= 1$ and $I= 1 \times {\rm SO}(Q^\perp)  \simeq {\rm SO}(6)$. \ps

As the subgroup $I$ fixes pointwise $L_1$ by construction, 
we have $I \subset \rho(T_1) = {\rm SO}(L_1^{\perp} )$, 
and thus its inverse image $I'=\rho^{-1}(I)$ in ${\rm Spin}(8)$ 
is included in $T_1 \simeq {\rm Spin}(7)$. 
If $M \subset T_1$ denotes the stabilizer of $L_3$, then $M$ is a ${\rm Spin}(1,6)$-subgroup
of $T_1$, and we have $I'=M$ in case (b1) and $I'=M^0$ in case (b2).
We also have $$\widetilde{S}_1 \cap \widetilde{S_2} = \rho^{-1}(\{ \pm 1\} \times I) = Z I'\,\, \, \, {\rm and}\, \, \,  S_1 \cap S_2 = \pi(\widetilde{S}_1 \cap \widetilde{S_2}),$$ 
so $S_1 \cap S_2 = \pi(Z I') = \pm \pi(I')  = \pi(I')$ as $\pi(\ker \rho)=\{ \pm 1\}$ and $\ker \rho \subset I'$.
It follows that $\pi$ induces isomorphisms $T_1 \isomo S_1$ and $I' \isomo S_1 \cap S_2$.
This shows that assertions (ii) and (iii) of the statement hold respectively in sub-cases (b1) and (b2).
\end{pf}
\begin{prop} 
\label{prop:accG2}
Let $r : \Gamma \rightarrow {\rm Spin}(7)$ be a morphism such that $r(\Gamma)$ 
is included in $\{ \pm 1\} \times H$ with $H \simeq {\rm G}_2$. Then $r$ is acceptable.
\end{prop}

\begin{pf}
We may choose $\eta$ such that $(r,\eta)$ satisfies (U1) and (U2). 
As the spin representation of ${\rm Spin}(7)$ is real, we may and do view it as an $8$-dimensional
$\R[{\rm Spin}(7)]$-module.
As $\R[H]$-modules, 
the representations $S$ and $E$ of ${\rm Spin}(7)$ satisfy $S \simeq E \oplus 1$.
As a $\{ \pm 1\} \times H$-module, we have thus $ S \,\simeq \,\eta \otimes E  \, \oplus\, \eta$, 
where $\eta$ denotes the first projection 
(whose restriction to $\Gamma$ is indeed the character already denoted by $\eta$),
and so the same isomorphism holds as $\R[\Gamma]$-modules.
But we also have $S \simeq \, \eta \otimes S$ by (U1) and Proposition~\ref{prop:propequat}, and thus 
$$E \oplus 1\, \simeq \, \eta \otimes E \, \oplus \,\eta.$$
As $\eta \neq 1$, and by semi-simplicity, this forces $\eta$ to appear as an $\R[\Gamma]$-submodule of $E$,
in contradiction with (U2) and Corollary~\ref{cor:coroetastable}.
\end{pf}

We are now able to prove Theorem~\ref{thm:exceptionPin6}.

\begin{pf} (of Theorem~\ref{thm:exceptionPin6})
Assume we have $r : \Gamma \rightarrow {\rm Spin}(7)$ and $\eta : \Gamma \rightarrow \{ \pm 1\}$ with $(r,\eta)$ unacceptable. Set $r' = \eta r$. We may choose a subgroup $S_1$ of ${\rm SO}(8)$ isomorphic to 
${\rm Spin}(7)$ and assume that $r$ is $S_1$-valued.  By (U1) and the acceptability of ${\rm O}(8)$, there is $g \in {\rm O}(8)$ such that $g r' g^{-1} = r$. 
In particular, we have 
\begin{equation} 
\label{eq:inclrgamma} 
r(\Gamma) \subset  S_1 \,\cap S_2 \,\, \, {\rm with}\,\, S_2 := g \,S_1\, g^{-1}. 
\end{equation}

Assume first $S_1=S_2$. Then $g$ is in the normalizer of ${\rm Spin}(7)$ in ${\rm SO}(8)$. 
But this normalizer is ${\rm Spin}(7)$,
as the latter only has inner automorphisms and centralizer $\{ \pm 1\}$ in ${\rm SO}(8)$. 
So we have $g \in {\rm Spin}(7)$, a contradiction as $(r,\eta)$ is unacceptable. 
This shows $S_1 \neq S_2$, and so we may apply Proposition~\ref{prop:subspin7}. It implies that 
inside $S_1 \simeq {\rm Spin}(7)$, the subgroup
$S_1 \cap S_2$ is either the normalizer of a ${\rm G}_2$-subgroup
or is included in a ${\rm Spin}(1,6)$-subgroup.
We conclude as the first case is excluded by Proposition~\ref{prop:accG2}.
\end{pf}

By Theorem~\ref{thm:exceptionPin6} and Proposition~\ref{prop:Nfus}, 
it is equivalent to classify the unacceptable $N$-valued or ${\rm Spin}(7)$-valued morphisms.

\begin{definition} 
\label{def:defYr} 
For any morphism $r : \Gamma \rightarrow {\rm Spin}(n)$ we denote by ${\rm Y}(r)$
the subset of ${\rm X}(r)$ consisting of characters of $\Gamma$ on the $1$-dimensional $\R[\Gamma]$-submodules of $E$.
\end{definition}

Another formulation of Theorem~\ref{thm:exceptionPin6} is:

\begin{coro} 
\label{cor:yrnonempty}
If $r : \Gamma \rightarrow {\rm Spin}(7)$ is unacceptable, we have ${\rm Y}(r)\neq \emptyset$.
\end{coro}

\section{Type I unacceptable morphisms}
\label{sect:typeI}

\begin{definition}
\label{def:typeI} 
Let $r : \Gamma \rightarrow {\rm Spin}(7)$ be unacceptable. 
We say that $r$ is {\rm of type I} if we have $1 \in {\rm Y}(r)$, or equivalently,
if the $\R[\Gamma]$-module $E$ contains $1$.
\end{definition}

%
%

Our goal in this section is to classify the type {\rm I} unacceptable ${\rm Spin}(7)$-valued morphisms.
We thus focus on the morphisms $r : \Gamma \rightarrow N$ satisfying 
$r(\Gamma) \subset N^0={\rm SU}(4)$, or equivalently, on the complex $4$-dimensional unitary representations of $\Gamma$ with determinant $1$.
We start with an important example, inspired by the analysis in \cite[\S 4.6]{chg2}, 
and that will turn out to be universal.

\ps
\begin{center} {\sc Example 2}\end{center}
\ps

Fix first a complex, non degenerate, quadratic plane $P \simeq \C^2$ 
and consider its similitude orthogonal group ${\rm GO}(P)$. 
By the choice of an orthonormal basis of $P$ we may and do identify $P$
with $\C^2$ equipped with the quadratic form $(x,y) \mapsto x^2+y^2$, in which case 
${\rm GO}(P)$ is the {\it standard} ${\rm GO}_2(\C)$.
Let $\mu : {\rm GO}(P) \rightarrow \C^\times$ be the similitude factor 
and $\det$ the determinant. The two isotropic lines of $P$ are 
permuted by ${\rm GO}(P)$, and we let $\epsilon : {\rm GO}(P) \rightarrow \{\pm 1\}$ 
be the signature of this permutation representation; its kernel is the subgroup 
${\rm GSO}(P) \simeq \C^\times \times \C^\times$ of proper similitudes. 
The structure of ${\rm GO}(P)$ is clear, 
for instance we have ${\rm ker}\, \mu\, = \,{\rm O}(P)$, 
$\det \, = \,\mu \, \epsilon$ and 
${\rm GO}(P)\, =\, \C^\times \cdot {\rm O}(P)$. 
We denote by ${\rm GO}(2)$ the (actually unique) 
maximal compact subgroup of ${\rm GO}_2(\C)$,
and define
$${\rm O}(2)^{\pm} \subset {\rm GO}(2) \subset {\rm GO}_2(\C)$$
as the subgroup of elements $g \in {\rm GO}(2)$ with $\mu(g) = \pm 1$.
The group ${\rm O}(2)^{\pm}$ is generated by the homothety $i {\rm 1}_2$ and ${\rm O}(2)$, 
hence is also isomorphic to the quotient of $\mu_4 \times {\rm O}(2)$ by its diagonal central $\Z/2\Z$. 
It has exactly $3$ order $2$ characters, 
namely $\det, \mu$ and $\epsilon$. 
Its tautological $2$-dimensional representation $P$ satisfies 
$\det P = \det$ and $P^\ast \simeq P \otimes \mu$ 
(with symmetric pairing); it is irreducible, non self-dual, 
and satisfies $P \simeq P \otimes \epsilon$.  \ps\ps

\begin{definition} 
\label{def:calG}
We denote by ${\mathcal{G}}$ the subgroup of all elements $(g_1,g_2)$ 
in ${\rm O}(2)^{\pm} \times {\rm O}(2)^{\pm}$ such 
that $\det g_1 = \det g_2$ and $\mu(g_2) = \epsilon(g_1)$. 
\end{definition}

The group ${\mathcal{G}}$ has two natural $2$-dimensional 
complex representations $P_1$ and $P_2$, given by the tautological 
representations of the first and second factors respectively. 
Fix an embedding $\rho : {\mathcal{G}} \rightarrow {\rm SU}(4)$ 
with underlying $4$-dimensional representation 
$\simeq P_1 \oplus P_2$. Such a $\rho$ is unique up to conjugacy by Lemma~\ref{lem:straccclassical},
but to fix ideas we take $P_1=\C^2 \times 0$ and $P_2 = 0 \times \C^2$
with standard quadratic and hermitian forms on each factor (namely $(x,y) \mapsto x^2+y^2$ and 
$|x|^2+|y|^2$).
We define morphisms $\upsilon$ and ${\rm d}: {\mathcal{G}} \rightarrow \{ \pm 1\}$ by setting, 
for $g=(g_1,g_2) \in \mathcal{G}$,
 \begin{equation}\label{def:upsd} \upsilon(g) = \epsilon(g_1)=\mu(g_2)\, \, \, {\rm and}\, \, \, {\rm d}(g) =\det(g_1)=\det(g_2).\end{equation}

\begin{prop} 
\label{prop:ex2Gamma}
View $\rho$ as a morphism ${\mathcal{G}} \rightarrow N$. Then $\rho$ and $\upsilon \rho$ 
are element conjugate in $N$, but not $N$-conjugate.
\end{prop}

\begin{pf}
By the discussion above, we have $P_1 \otimes \upsilon \simeq P_1$ 
and $P_2 \otimes \upsilon \simeq P_2^\ast$, hence
\begin{equation}\label{propchiphi} \rho \simeq P_1 \oplus P_2 \hspace{1 cm} {\rm and}\hspace{1 cm} \upsilon \rho \simeq \upsilon \otimes \rho \simeq P_1 \oplus P_2^\ast. \end{equation} 
As neither $P_1$ nor $P_2$ is self-dual, neither $\rho$ nor its dual is isomorphic to $ \upsilon \rho$, so $\rho$ and $\upsilon \rho$ are not $N$-conjugate, by Remark~\ref{rem:conjinsN}. However, for all $g \in {\mathcal{G}}$ we have $\rho(g)$ or $\rho(g)^{-1}$ conjugate to $\upsilon(g) \rho(g)$ in ${\rm SU}(4)$. 
Indeed, this is trivial for $g$ in ${\rm ker}\, \upsilon$. 
It is thus enough to show that on both ${\rm ker}\, {\rm d}$ 
and ${\rm ker}\, \upsilon {\rm d}$ we have $\rho^\ast \simeq \upsilon \rho$, and for that it is enough to prove that on those two subgroups we have $P_1 \simeq P_1^\ast$. But on ${\rm ker}\, {\rm d}$, this follows from the fact that we have $1={\rm d}(g)= \det(g_1)$ ($P_1$ is symplectic), and on ${\rm ker}\, \upsilon {\rm d}$, this follows from $1=({\rm d} \upsilon)(g) = \mu(g_1)$ ($P_1$ is orthogonal). 
\end{pf}
\ps\ps

We now start showing that any unacceptable morphism $r : \Gamma \rightarrow N$ with $r(\Gamma) \subset N^0$
can be explained by this example. We denote by $V$ the tautological complex $4$-dimensional representation of $N^0={\rm SU}(4)$. For a given morphism $r : \Gamma \rightarrow N$ with $r(\Gamma) \subset N^0$, we may view $V$ as a (semi-simple) $\C[\Gamma]$-module.

\begin{lemma}
\label{lem:candecV}
Let $r : \Gamma \rightarrow {\rm SU}(4)$ and $\eta : \Gamma \rightarrow \{ \pm 1\}$ be two morphisms. Assume $\eta\, r$ and $r$ are element conjugate in $N$
but not $N$-conjugate. Then the $\C[\Gamma]$-module $V$ defined by $r$ has a unique decomposition
$V = A \oplus B$ such that:
\begin{itemize}
\item[(i)] $\dim A=\dim B=2$,\ps
\item[(ii)] we have $\C[\Gamma]$-module isomorphisms $A \simeq \eta \otimes A$ and $B^\ast \simeq \eta \otimes B$.\ps
\end{itemize}
Moreover, the following properties hold:
\begin{itemize}
\item[(a)] neither $A$ nor $B$ is self-dual,\ps
\item[(b)] the character $\det A=\det B$ has order $2$, and we have $\det A \neq \eta$,\ps
\item[(c)] the $\C[\Gamma]$-module $V$ is multiplicity free,\ps
\item[(d)] $A$ is reducible if, and only if, we have $A \simeq a \oplus \eta a$ with $a$ of order $4$,\ps
\item[(e)] $B$ is reducible if, and only if, we have $B \simeq b_1 \oplus b_2$ with $b_1^2 = b_2^2 =\eta$.
\end{itemize}
\end{lemma}

\begin{pf} 
We start as in \cite[Lemma 2.4]{larsen2}.
Denote respectively by $V_1$ (resp. $V_2$) 
the $4$-dimensional representation of $\Gamma$ on $V$ defined by $r$ (resp. $\eta r$). 
We have thus 
\begin{equation}\label{eq:v2v1eta} V_2 \simeq \eta \otimes V_1. \end{equation}
and $\det V_i=1$. The $V_i$ are semi-simple.
By (U1) and Remark~\ref{rem:conjinsN} we have
$V_1 \oplus V_1^\ast \simeq V_2 \oplus V_2^\ast$. 
We may write  
$V_1 \simeq  A \oplus B$ and $V_2 \simeq A \oplus C$ with ${\rm Hom}_\Gamma(B,C)=0$. 
The previous relation shows $B \oplus B^\ast \simeq C \oplus C^\ast$ 
and then $C \simeq B^\ast$. 
In particular, we have 
$\det A\, \det B \,= \,\det A \,(\det B)^{-1}\, =1$, so 
$$c:=\det A = \det B$$
 is a character of $\Gamma$ of order $\leq 2$. Also, 
(U2) is equivalent to $V_1 \not \simeq V_2$ and $V_1 \not \simeq V_2^*$, {\it i.e.} 
to: neither $A$ nor $B$ is self-dual. This proves property (a) of the statement.
In particular, both $A$ and $B$ have dimension $\geq 2$, 
hence must have dimension $2$. As a summary, we have shown 
\begin{equation}\label{eq1} V_1 \simeq A \oplus B, \, V_2 \simeq A \oplus B^\ast, \, A^\ast \not \simeq A, \, B^\ast \not \simeq B, \, \, \dim A = \dim B = 2.\end{equation}
In particular, the character $c$ must be non trivial, since any $2$-dimensional representation $P$ satisfies $P^\ast \otimes \det P \simeq P $. So $c$ has order $2$.  \ps
By \eqref{eq:v2v1eta} we have $A \oplus B^\ast \simeq \eta \otimes (A \oplus B)$.
Note that if two semi-simple $2$-dimensional representations 
of same determinant share a $1$-dimensional constituent, 
then they are isomorphic. As a consequence, 
if $A$ is not isomorphic to $\eta \otimes A$, 
then we have $A \simeq \eta \otimes B$ 
and $B^\ast \simeq \eta \otimes A$. But these two relations imply 
$A \simeq A^\ast$, a contradiction. We have proved that the decomposition 
$V=A \oplus B$ satisfies (i) and (ii).\ps
Note that $A \simeq \eta \otimes A$ is equivalent to 
$A^\ast \simeq c \eta \otimes A$. In particular, we have $c \neq \eta$ as $A$ is not selfdual,
and this ends proving property (b).
If $A$ is reducible and contains say the character $a : \Gamma \rightarrow \C^\times$, it also contains $\eta a \neq a$, and
we have $A \simeq a \oplus \eta a$ with $c= a^2 \eta$, $a^2 \neq 1$ and $a^4=c^2=1$, hence property (d).
Similarly, if $B$ is reducible and contains $b$, we have $B \simeq b \oplus c b^{-1}$
and the relations $B^\ast \simeq \eta \otimes B$ and $c \neq \eta$ imply $b^{-1} = \eta b$.
This proves property (e).
In any case, both $A$ and $B$ are multiplicity free as we have $\eta, c \eta \neq 1$.
Also, there is no nonzero 
$\Gamma$-equivariant morphism $A \rightarrow B$, since such a morphism would be an isomorphism and would imply $B \simeq \eta \otimes B \simeq B^\ast$.
We have proved property (c): $V$ is multiplicity free.
\ps
It only remains to show the uniqueness statement. 
Assume we have a $\C[\Gamma]$-module decomposition $V= A' \oplus B'$
with $A'$ and $B'$ satisfying (i) and (ii). 
Consider the $\Gamma$-equivariant projection $f: A' \rightarrow B$ with kernel $A$.
It cannot be an isomorphism,
since it would imply $B \simeq \eta \otimes B$ by (ii) and then $B^\ast \simeq B$.
If $f$ is nonzero, then both $A$ and $B$ are reducible and we have 
$A' \simeq a \oplus b$ with $a$ (resp. $b$) a constituent of $A$ (resp. $B$),
and in particular $a^2 \eta = c$ and $b^2 = \eta$. But $A' \simeq \eta \otimes A'$
implies $b=\eta a$, so we have $c=\eta a^2 = \eta b^2 = 1$: a contradiction.
We deduce $f(A')=0$, {\it i.e.} $A' \subset A$, and then $A'=A$ and $B' \simeq B$.
As $V$ is multiplicity free, this forces $B'=B$ as well.
\end{pf}

\begin{remark}
If both $A$ and $B$ are reducible, it follows from Lemma~\ref{lem:candecV} that we have
$A \simeq a \oplus a b^2$, $B  \simeq b \oplus b a^2$, $a,b$ of order $4$ and $\eta = b^2$.
This is thus exactly the case studied in {\sc Example 1}.
\end{remark}

For later use we mention the following complement to Lemma~\ref{lem:candecV}.

\begin{lemma}
\label{compl:l2biggamma}
In the notations of Lemma~\ref{lem:candecV}, the $\C[\Gamma]$-module  $A \otimes B$ does not contain any character $\chi: \Gamma \rightarrow \C^\times$ with $\chi^2=1$.
\end{lemma}

\begin{pf}
Indeed, such a $\chi$ induces a nonzero $\C[\Gamma]$-equivariant morphism 
$\chi \otimes A^\ast \rightarrow B$. 
This forces $\chi \otimes A^\ast \simeq B$ 
since both $\chi \otimes A^\ast$ and $B$ have determinant $\chi^2 c^{-1} = c$, as $\chi^2=1$. 
This implies $c \chi \otimes A \simeq B$ and then ${\rm Sym}^2 A \simeq {\rm Sym}^2 B$.
But ${\rm Sym}^2 B$ contains $\eta$ by Lemma~\ref{lem:D8a} (i), 
hence so does ${\rm Sym}^2 A$. This forces
$\eta \otimes A^\ast \simeq A$ (again since both sides have the same determinant) and then $A^\ast \simeq A$: a contradiction.
\end{pf}

The following corollary applies for instance to $\Gamma = \Z \times {\rm SU}(2)$.

\begin{coro} 
\label{coro:homz22}
Assume there is no surjective morphism $\Gamma \rightarrow \Z/2\Z \times \Z/2\Z$.
Then any morphism $\Gamma \rightarrow {\rm Spin}(7)$ is acceptable.
\end{coro}

\begin{pf} 
If $(r,\eta)$ is unacceptable, the assumption implies ${\rm Hom}(\Gamma,\{ \pm 1\}) = \{1,\eta\}$.
But then ${\rm Y}(r)=\{1\}$ by Corollary~\ref{cor:coroetastable} and Theorem~\ref{thm:exceptionPin6}, 
so $r$ is of type I.
But this contradicts the property $\det A \neq 1,\eta$ in Lemma~\ref{lem:candecV}
\end{pf}

\ps\ps

\begin{thm}
\label{thm:univGG}
Assume $r,r' : \Gamma \rightarrow N$ are element conjugate but non conjugate morphisms
with values in ${\rm SU}(4)=N^0$. Then there exists a morphism 
$f : \Gamma \rightarrow {\mathcal{G}}$ such that, 
up to replacing $r$ and $r'$ by some $N$-conjugate if necessary, we have: \ps 
\begin{itemize}
\item[(i)] $r= \rho \circ f $, \ps
\item[(ii)] $r'=\eta \,r $ with $\eta:=\upsilon \circ f$, and\ps
\item[(iii)] $\eta$ and ${\rm d} \circ f$ are distinct, order $2$, characters of $\Gamma$.
\end{itemize}
\end{thm}

\begin{pf} 
By Proposition-Definition~\ref{def:defetaN}, up to conjugating $r'$ by some element of $N$ 
we may assume that we have $r' = \eta r$ 
for some order $2$ character $\eta : \Gamma \rightarrow \{\pm 1\}$. 
Write $V \simeq A \oplus B$ as in Lemma~\ref{lem:candecV} and set $c=\det A$.
Apply Lemma~\ref{lem:D8a} (i) below to both $(P,\eta)=(A,c\eta)$ and $(P,\eta)=(B,\eta)$.  
It endows $A$ and $B$ with nondegenerate $\Gamma$-equivariant 
symmetric pairings ${\rm b}_A$ and ${\rm b}_B$ 
with respective similitude factors $c\eta$ and $\eta$. 
We may assume that both quadratic spaces $(A,{\rm b}_A)$ and $(B,{\rm b}_B)$ are 
$\C^2$ equipped with the standard form $x^2+y^2$. 
The action of $\Gamma$ on $A$ and $B$ thus gives rise
to morphisms $f_A$ and $f_B: \Gamma \rightarrow {\rm GO}_2(\C)$,
with respective similitude factors $\mu \circ f_A = c \eta$ and $\mu \circ f_B =\eta$,
and $\det \circ f_A = \det \circ f_B = c$.
But $\Gamma$ also preserves some positive definite hermitian form inherited from $V$, 
hence $f_A(\Gamma)$ and $f_B(\Gamma)$ have a compact closure in ${\rm GO}_2(\C)$, 
and so are included in ${\rm GO}(2)$, and even in ${\rm O}(2)^{\pm}$ as we have 
$\eta^2=(c \eta)^2=1$. We have thus constructed a morphism $f := f_A \times f_B : \Gamma \rightarrow {\mathcal{G}}$ such that 
the two morphisms $\rho \circ f$ and $r$, from $\Gamma$ to ${\rm SU}(4)$, define two 
isomorphic representations of $\Gamma$ on $\C^4$. 
By Lemma~\ref{lem:straccclassical}, they are conjugate in ${\rm SU}(4)$.
So up to replacing $r$ with a conjugate 
we may assume $\rho \circ f = r$. The result follows from the formulae 
$\upsilon \circ \rho = \mu \circ f_B = \eta$ and ${\rm d} \circ \rho = \det \circ f_A = c$.
\end{pf}

In the following lemma, we denote by ${\rm D}_8$ the dihedral group of order $8$.

\ps
\begin{lemma}
\label{lem:D8a}\label{lem:D8b}
 Let $\Gamma$ be a group, $P$ a $2$-dimensional complex representation of $\Gamma$ and $\eta : \Gamma \rightarrow \C^\times$ a morphism. Assume $\det P \neq \eta$ and $P^\ast  \simeq P \otimes \eta^{-1}$. \ps
 \begin{itemize}
 \item[(i)] There is a nondegenerate symmetric pairing on $P$ such that $\Gamma$ acts on $P$ as orthogonal similitudes with similitude factor $\eta$, and ${\rm Sym}^2 P$ contains $\eta$.\ps
 \item[(ii)] Assume furthermore that $\det P$ and $\eta$ are two order $2$ characters of $\Gamma$. 
 Then $P$ is self-dual if, and only if, the image of $\Gamma$ in ${\rm GL}(P)$ is isomorphic to ${\rm D}_8$.
 \end{itemize}
 \end{lemma}

 \begin{pf} The given isomorphism $P^\ast  \simeq P \otimes \eta^{-1}$ may be viewed a nondegenerate $\Gamma$-equivariant pairing $P \otimes P \rightarrow \eta$. As we have $\Lambda^2 P = \det P \, \neq \eta$ by assumption, this pairing factors through ${\rm Sym}^2 P$. This proves (i). We now prove assertion (ii). \par
 Of course, the unique faithful $2$-dimensional representation of ${\rm D}_8$ is self-dual. 
Conversely, assume $P^\ast \simeq P$. We have then
$P^\ast \simeq P \otimes \nu$ for each of the four distinct characters $\nu = 1, \eta, \det, \eta \det$, and thus ${\rm Sym}^2 P \simeq 1 \oplus \eta \oplus \eta \det$. Note that $(\eta,\det) : \Gamma \rightarrow \mu_2^2$ is surjective, and that the kernel of the natural morphism ${\rm GL}(P) \rightarrow {\rm GL}({\rm Sym}^2 P)$ is the central $\mu_2$. As a consequence, the image $I$ of $\Gamma$ in ${\rm GL}(P)$ is an extension of $\mu_2^2$ by a central subgroup of order $\leq 2$. This shows $|I|\leq 8$. Note also that $P$ cannot contain any character $\nu$, otherwise it would contain the three distinct characters $\nu$, $\nu \eta$, $\nu \det$. So $P$ is irreducible and $I$ is nonabelian of order $8$. Since the unique $2$-dimensional faithful representation of the quaternion group of order $8$ has determinant $1$, $I$ must be isomorphic to ${\rm D}_8$.
 \end{pf}
\ps\ps

We now state a converse to Theorem~\ref{thm:univGG}. 
\ps\ps

\begin{prop}
\label{prop:converseN}
Let $\Gamma$ be a group and $f : \Gamma \rightarrow {\mathcal{G}}$ a morphism
such that ${\rm d} \circ f$ and $\eta:=\upsilon \circ f$ are distinct and nontrivial. 
Then the $N$-valued morphisms $r:=\rho \circ f$ and $\eta r$ 
are element-conjugate by Proposition~\ref{prop:ex2Gamma}. 
They are non conjugate if, and only if, 
none of the two projections $\Gamma \rightarrow {\rm O}(2)^{\pm}$ 
has an image isomorphic to ${\rm D}_8$.
\end{prop}

\begin{pf} 
This follows from Remark~\ref{rem:conjinsN}, Formula \eqref{propchiphi} and Lemma~\ref{lem:D8b} (ii). 
\end{pf}

We finally provide an alternative study of the unacceptable morphisms of type I from the point of view of their
standard representation on $E$. 
View $\rho$ as a morphism ${\mathcal{G}} \rightarrow {\rm Spin}(7)$, 
and consider the associated $\R[{\mathcal{G}}]$-module $E$. We have 
a $\mathcal{G}$-stable decomposition $E = L \perp F$,
with $L \simeq 1$ as an $\R[\mathcal{G}]$-module, and 
by Remark~\ref{rem:repNE} we have a $\C[\mathcal{G}]$-module decomposition
\begin{equation}
\label{eq:decCGF}
F \otimes \C\,  \simeq \,\, P_1 \otimes P_2\,\, \oplus\, \,{\rm d} \,\oplus \, {\rm d}.
\end{equation}

\begin{prop}
\label{prop:descorthG}
We have an orthogonal decomposition 
\begin{equation}
\label{eq:decFj} 
F = F_1 \perp F_2 \perp F_3, \, \, \, {\rm with}\, \, \dim F_1 = \dim F_2 = \dim F_3 = 2,
\end{equation}
such that ${\mathcal{G}}$ coincides with the 
subgroup of elements $g \in N$ with determinant $1$ on $F$, 
preserving the pair $\{ F_1, F_2\}$ (hence $F_3$), 
and acting on $F_3$ as a homothety whose sign coincides with the signature of $g$ 
on the $2$-elements set $\{ F_1,F_2\}$. Moreover, this latter signature is ${\rm d}(g)$, and if 
we assume ${\rm d}(g)=1$ then we have $\upsilon(g) = \det_{F_1}(g) = \det_{F_2}(g)$. 
\end{prop}

In the decomposition above, both $F_3$ and the pair $\{F_1,F_2\}$ are thus unique, 
but of course the numberings of $F_1$ and $F_2$ are not.

\begin{pf} 
Consider the subgroup $T:= {\rm SO}(2) \times {\rm SO}(2)$ of $\mathcal{G}$ (a torus).
By Formula~\eqref{eq:decCGF},  
there is a unique triple $\{F_1,F_2,F_3\}$ of planes of $F$ with $F = F_1 \perp F_2 \perp F_3$, 
such that $T$ acts trivially on $F_3$ and the natural map $T \rightarrow {\rm SO}(F_1) \times {\rm SO}(F_2)$ is surjective, with kernel $\{\pm 1\}$. 
More precisely, if we write $P_1 \simeq x \oplus x^{-1}$ and $P_2 \simeq y \oplus y^{-1}$ as $\C[T]$-modules, then up to replacing $x$ with $x^{-1}$ we may assume $F_1\otimes \C \simeq xy \oplus (xy)^{-1}$ and $F_2 \otimes \C \simeq xy^{-1} \oplus x^{-1}y$. 
By Definition~\ref{def:calG}, the group ${\mathcal{G}}$ is generated by $T$ and
the two elements $g_1=(\sigma, i {\rm 1}_2)$ and $g_2=(i {\rm 1}_2,\sigma)$, for any $\sigma \in {\rm O}(2)$ with determinant $-1$. Both $g_1$ and $g_2$ exchange thus $F_1$ and $F_2$.
We conclude as we have $\upsilon(g_1)=-1$, $\upsilon(g_2)=1$, ${\rm d}(g_1)={\rm d}(g_2)=-1$, 
and $\det_{F_j}(g_1 g_2)=-1$ for $j=1,2$.
\end{pf}

Let ${\mathcal{G}}_1 \subset {\mathcal{G}}$ be the kernel of ${\rm d}$ and consider the $4$-dimensional Euclidean space
$$F'= F_1 \perp F_2.$$
By Proposition~\ref{prop:descorthG}, we may identify $F'$ with the natural direct sum $\R^2 \perp \R^2$, 
${\mathcal{G}}/\{ \pm 1\}$ with the subgroup $\underline{{\mathcal{G}}}$ of ${\rm SO}(4)$ 
preserving this sum,  
and ${\mathcal{G}_1}/\{ \pm 1\}$ with 
$$\underline{{\mathcal{G}_1}} :=  \{ (g_1,g_2) \in {\rm O}(2) \times {\rm O}(2)\, \, |\, \, \det(g_1)=\det(g_2)\}.$$ 
The element $\theta = {\tiny \left[ \begin{array}{cc} 0 & {\rm 1}_2 \\ {\rm 1}_2 & 0 \end{array}\right]}$ has determinant $1$ and is in 
$\underline{{\mathcal{G}}}\smallsetminus \underline{{\mathcal{G}_1}}$, with $\theta^2=1$. 
For all $(g_1,g_2) \in \underline{{\mathcal{G}_1}}$, we have $\theta(g_1,g_2) \theta^{-1} = (g_2,g_1)$ and $\upsilon(g_1,g_2)=\det(g_1)=\det(g_2)$ (the value of $\upsilon$ on $\theta$ will be irrelevant). \ps

We have proved that up to conjugacy, any type I unacceptable morphism $(r,\eta)$ arises 
from a morphism $\Gamma \rightarrow {\mathcal{G}}$ in the sense of Theorem~\ref{thm:univGG}. 
A first converse statement was given in Proposition~\ref{prop:converseN}, in terms of the 
natural morphism $\Gamma \rightarrow N^0 ={\rm SU}(4)$.
Our aim now is to provide a second one, rather in terms of the properties of the $\R[\Gamma]$-module $E$. 
We first highlight a few examples.

\begin{prop}
\label{exa:nonexamples}
We denote by $\Delta \simeq {\rm O}(2)$ the diagonal subgroup of ${\rm O}(2) \times {\rm O}(2)$ 
and by $\mu \simeq \{ \pm 1\}^4$ the diagonal subgroup of ${\rm O}(4)$. 
We define three subgroups $H_1,H_2$ and $H_3$ of ${\rm SO}(4)$ 
contained in $\underline{\mathcal{G}}$ as follows:\par
\begin{itemize}
\item[(i)] $H_1$ is generated by $\mu \cap \Delta$ and ${\rm diag}(s,-s) \theta$ with 
$s={\tiny \left[ \begin{array}{cc}  0  & 1 \\ 1 & 0 \end{array}\right]}$, \par
\item[(ii)] $H_2$ is generated by $\Delta$ and $\theta$,\ps
\item[(iii)] $H_3$ is generated by $\Delta$ and $t:={\rm diag}({\rm 1}_2,-{\rm 1}_2) \theta$.\ps
\end{itemize}
\par
For each $i$, the characters $\delta = {\rm d}_{|H_i}$ and $\eta={\rm \upsilon}_{|H_i}$ of $H_i$ are distinct and of order $2$, and we denote by $r$ the inclusion of $H_i$ in ${\rm Spin}(7)$. 
Then $(r,\eta)$ is acceptable for $i=1,2$ and unacceptable for $i=3$.
\end{prop}

\begin{pf}
Each $H_i$ contains elements of the form $g$ and $g'\theta$ with $g,g' \in \Delta$ and $\upsilon(g)=-1$, 
hence the first assertion.
We have $H_1 \simeq {\rm D}_8$ and $H_2 \simeq {\rm O}(2) \times \Z/2\Z$, and for $i=1,2$, 
 an $\R[H_i]$-module isomorphism $F'\, \simeq\,\, \delta \otimes P \,\,\oplus \,\,P$ 
with $P$ irreducible and $\det P \in \{ \eta, \eta \delta\}$. So $\eta \in {\rm X}(r)$ and $(r,\eta)$ is acceptable for $i=1,2$. 
On the other hand, we have $H_3 \simeq ({\rm O}(2) \times \mu_4)/{\rm diag} \,\mu_2$ and an 
$\R[H_3]$-module isomorphism $F' \simeq P \otimes Q$, with $P$ the (inflation of the) tautological representation of the ${\rm O}(2)$ factor, and $Q$ (that of) the irreducible 2-dimensional real representation of the $\mu_4$ factor. Indeed, the element $t$ commutes with $\Delta$ and satisfies $t^2={\rm diag}(-{\rm 1}_2,-{\rm 1}_2)\in \Delta$. So $F'$ is irreducible, we have 
${\rm X}(r)=\{1,\delta\}$, and $(r,\eta)$ is unacceptable for $i=3$.
\end{pf}

The two first examples turn out to explain all the acceptable examples.
For minor reasons we first need to introduce the subgroup 
$M \subset N$ preserving $F_1$, $F_2$ (hence $F_3$) and acting trivially on $F_3$. We have 
$M/\{ \pm 1\} \simeq {\rm O}(2) \times {\rm O}(2)$, $M \cap {\mathcal{G}} = {\mathcal{G}}_1$ and
$M$ normalizes ${\mathcal{G}}$ by fixing the character ${\rm d}$ and exchanging $\upsilon$ and $\upsilon {\rm d}$. 

\begin{prop} 
\label{prop:classnonex}
Consider a morphism $f : \Gamma \rightarrow {\mathcal{G}}$ such that 
$\eta = \upsilon \circ f$ and $\delta = {\rm d} \circ f$ are distinct and of order $2$,
and set $r = \rho \circ f$. Then $(r,\eta)$ satisfies {\rm (U1)}. It is acceptable if, and only if, 
up to conjugating $f$ with some element of $M$, 
the image of $f(\Gamma)$ in $\underline{{\mathcal{G}}}$ is included in one of the subgroups
$H_1$ and $H_2$ of Proposition~\ref{exa:nonexamples}.
\end{prop}

\begin{pf} 
That $(r,\eta)$ satisfies (U1) follows from Proposition~\ref{prop:ex2Gamma},
so we only have to check the last assertion. 
The sufficiency of the given condition follows directly from
Proposition~\ref{exa:nonexamples}, so we assume from now on that $(r,\eta)$ is acceptable, 
or equivalently, that we have $\eta \in {\rm X}(r)$. Set $\Gamma_1 = \ker \delta$. \ps

By~\eqref{eq:decFj} we may find some $\R[\Gamma]$-submodule $U \subset F'$
with $\det_U \in \{ \eta,\eta \delta\}$. Il all cases $\det_U$ coincides with $\eta$ on $\Gamma_1$,
and up to replacing $U$ by its orthogonal in $F'$ if necessary, we may also assume $1 \leq \dim U \leq 2$. 
Note that $U$ is not included in $F_1$ or in $F_2$,
since $\Gamma$ permutes transitively $\{F_1,F_2\}$ as $\delta \neq 1$.
So the natural $\Gamma_1$-equivariant projection $U \rightarrow F_i$ has a nonzero image.
As $U,F_1$ and $F_2$ have the same determinant $\eta$ on $\Gamma_1$,
this implies that $F_1$ and $F_2$ are isomorphic as $\R[\Gamma_1]$-modules (and 
isomorphic to $U_{|\Gamma_1}$ in the case $\dim U=2$). 
By Lemma~\ref{lem:straccclassical}, and up to replacing $f$ by some $M$-conjugate, we may thus assume $f(\Gamma_1) \subset \Delta$. (We freely use the notations of Proposition~\ref{exa:nonexamples}). \ps

As $\eta$ is nontrivial on $\Gamma_1$, 
there is an element $h \in \Gamma$ with $f(h)={\rm diag}(a,b)\theta$, $a,b \in {\rm O}(2)$ 
and $\det a = \det b = -1$. 
The element $f(h^2)=f(h)^2 = {\rm diag}(ab,ba)$ is in $\Delta$, so the two reflections 
$a$ and $b$ in ${\rm O}(2)$ commute. This forces $b=\epsilon a$ for some sign 
$\epsilon = \pm 1$. If $\epsilon=1$ then we have $f(\Gamma)\subset H_2$ and we are done.
So we definitely assume $\epsilon = -1$, {\it i.e.} $f(\Gamma) \subset H_3$, and we 
write $H_3 \simeq ({\rm O}(2) \times \mu_4)/{\rm diag} \,\mu_2$ and $F' \simeq P \otimes Q$ as in the proof of Proposition~\ref{exa:nonexamples}. \ps
Assume first $P$ is irreducible as an $\R[\Gamma_1]$-module (necessarily absolutely irreducible,
since $\det P \neq 1$). We claim that $F'$ is an irreducible $\R[\Gamma]$-module.
Indeed, we have $(F')_{|\Gamma_1} \simeq P \otimes Q$ with $Q$ trivial, 
so the proper $\R[\Gamma_1]$-submodules of $F'$  are the $P \otimes v$ with $v \in Q$ nonzero. 
But none of those is stable by the element $f(h) \in \Delta t$, since $t$ has no eigenvector in $Q$.
So we have $\det F'=1$, ${\rm X}(r) = \{1,\delta\}$ and $\eta \notin {\rm X}(r)$, contradicting the acceptability of $(r,\eta)$.\ps

So the $\R[\Gamma_1]$-module $P$ is reducible, necessarily the sum of two distinct real characters 
as $\det P = \eta \neq 1$. Up to replacing $f$ by some $M$-conjugate, we may thus assume $f(\Gamma_1) \subset \mu \cap \Delta$. Since $f(h)={\rm diag}(a,-a)\theta$ normalizes $f(\Gamma_1)$, 
the reflection $a \in {\rm O}(2)$ normalizes the diagonal subgroup $\{ \pm 1\}^2$ of ${\rm O}(2)$, so we either 
have $a=\pm {\tiny \left[ \begin{array}{cc}  0  & 1 \\ 1 & 0 \end{array}\right]}$ or 
$a = \pm {\tiny \left[ \begin{array}{cc}  1  & 0 \\ 0 & -1 \end{array}\right]}$. In the first case we have 
$f(\Gamma) \subset H_1$ and we are done. In the second case, $P$ is reducible as $\R[\Gamma]$-module, say $P \simeq \alpha \oplus \beta$, and we have $F' \simeq \alpha \otimes Q \oplus \beta \otimes Q$, a sum of two irreducible 
$\R[\Gamma]$-modules with determinant $1$.  This shows ${\rm X}(r)=\{1,\delta\}$ and contradicts
the acceptability of $(r,\eta)$.
\end{pf}

We end this paragraph by discussing a few other properties of type {\rm I} unacceptable morphisms.

\begin{coro} 
\label{cor:uniquedelta}
Assume $r : \Gamma \rightarrow {\rm Spin}(7)$ is unacceptable of type {\rm I}.  
We have ${\rm Y}(r)=\{1,\delta\}$ where $1$ and $\delta$ have respective 
multiplicity $1$ and $2$ in $E$. 
\end{coro}

\begin{pf} This is a corollary of Proposition~\ref{prop:Nfus}, Theorem~\ref{thm:univGG} and Lemma~\ref{compl:l2biggamma}.
\end{pf}

The next proposition asserts that for an unacceptable morphism $r$ of type {\rm I},
the action of ${\rm X}(r)$ on ${\rm E}(r)$ discussed in Remark~\ref{rem:etar} is transitive.

\begin{prop} 
\label{prop:etacan}
Assume $r : \Gamma \rightarrow {\rm Spin}(7)$ is unacceptable of type {\rm I}.
Then up to the multiplication by an element of ${\rm X}(r)$, there is a unique 
morphism $\eta : \Gamma \rightarrow \{ \pm 1\}$ such that $(r,\eta)$ 
is unacceptable.
\end{prop}

\begin{pf} Let $\delta$ be as in Corollary~\ref{cor:uniquedelta} and set $\Gamma_1=\ker \delta$.
Choose $\eta$ such that $(r,\eta)$ is unacceptable, and set $\eta_1=\eta_{|\Gamma_1}$.
By Theorem~\ref{thm:univGG} and Proposition~\ref{prop:descorthG} we have an $\R[\Gamma]$-module isomorphism $E \simeq 1 \oplus \delta \oplus \delta \oplus {\rm ind}_{\Gamma_1}^\Gamma U$,
where $U$ is an $\R[\Gamma_1]$-module of dimension $2$ with $\det U = \eta_1$. 
If $U$ is irreducible, then the $2$-dimensional irreducible constituents of the $\R[\Gamma_1]$-module
$E$ are $U$ and its outer conjugate under $\Gamma/\Gamma_1$, and both have determinant $\eta_1$.
It follows that $\eta_1$ is uniquely defined, hence that $\eta$ is unique up to multiplication by $\delta \in {\rm X}(r)$ in this case. So we may assume $U$ is reducible, say $U \,\simeq\, \alpha \oplus \alpha \eta_1$ for some character
$\alpha : \Gamma_1 \rightarrow \{ \pm 1\}$. 
We have $${\rm ind}_{\Gamma_1}^\Gamma U \simeq V \oplus \eta \otimes V, \, \, {\rm with}\,\,
V \,:=\, {\rm ind}_{\Gamma_1}^\Gamma \alpha.$$ 
\ps Denote by $\alpha^c$ the outer-conjugate of $\alpha$
under $\Gamma/\Gamma_1$. Let us assume first $\alpha^c \neq \alpha$. 
Then the image of $\Gamma_1$ and $\Gamma$ in ${\rm O}(V)$ 
are thus respectively isomorphic to $\mu_2^2$ and ${\rm D}_8$. 
The three order $2$ characters of this ${\rm D}_8$ must be $\delta$, $\delta \det_V$ and $\det_V$, 
and we have ${\rm X}(r)=\langle \delta, \det_V \rangle$.
As $r$ is unacceptable, $\eta$ is not in ${\rm X}(r)$, so the group 
$\pi(r(\Gamma))$ is isomorphic to ${\rm D}_8 \times \mu_2$. 
By Proposition~\ref{prop:facteta}, $\eta$ is an order $2$ character of this group.
We conclude as there is a unique non trivial character of $\pi(r(\Gamma))$ modulo ${\rm X}(r)$.
Assume now $\alpha^c=\alpha$. Then the image of $\Gamma_1$ 
and $\Gamma$ in ${\rm O}(V)$ are respectively isomorphic to $\mu_2$ and $\mu_4$. 
Indeed, that of $\Gamma$ cannot be $\mu_2 \times \mu_2$,
otherwise $V$ would contain a $1$-dimensional representation 
and contradict Corollary~\ref{cor:uniquedelta}. 
So we have ${\rm X}(r)=\{1,\delta\}$, 
and $\pi(r(\Gamma)) \simeq \mu_4 \times \mu_2$ has a unique
order $2$ character modulo ${\rm X}(r)$.
(This case is the one of {\sc Example 1}).
\end{pf}

\section{Type II unacceptable morphisms}
\label{sect:typeII}

\begin{definition} 
\label{def:deftypeII}
Let $r : \Gamma \rightarrow {\rm Spin}(7)$ be unacceptable, $\chi \in {\rm Y}(r)$
{\rm non trivial}, and set $\Gamma_0 = \ker \chi$ (an index $2$ subgroup of $\Gamma$). 
We say that $r$ is of type {\rm II} with respect to $\chi$ if there is $\eta : \Gamma \rightarrow \{\pm 1\}$ such that both $(r,\eta)$ and $(r_{|\Gamma_0},\eta_{|\Gamma_0})$ are unacceptable. 
In this situation, we also say that $(r,\eta)$ is unacceptable of type {\rm II} with respect to $\chi$.
We shall also say that $r$ is of type {\rm II} if it is so with respect to some $\chi$.
\end{definition}

Of course, to check the unacceptability of both $(r,\eta)$ and $(r_{|\Gamma_0},\eta_{|\Gamma_0})$, it is enough to check that $(r,\eta)$ satisfies (U1)
and that $(r_{|\Gamma_0},\eta_{|\Gamma_0})$ satisfies (U2). 
The following remark shows that types ${\rm I}$ and ${\rm II}$ are exclusive.

\begin{remark}
Assume that $r : \Gamma \rightarrow {\rm Spin}(7)$ is unacceptable of type {\rm I}.
By Corollary~\ref{cor:uniquedelta}, we have ${\rm Y}(r)=\{1,\delta\}$ with $1$ (resp. $\delta$) occuring with multiplicity $1$
(resp. $2$). So $r_{|\ker \delta}$ contains $1$ with multiplicity $2$, and $r$ is not of type {\rm II}.
\end{remark}
\ps

{\bf Notation} : Up to the end of this section, we assume that $\Gamma$ is a group, 
that $\chi : \Gamma \rightarrow \{ \pm 1\}$ is a non trivial character, and we 
denote by $\Gamma_0$ the kernel of $\chi$, an index $2$ subgroup of $\Gamma$. 
\ps


\begin{lemma} 
\label{lem:resGo}
Let $r : \Gamma \rightarrow {\rm Spin}(7)$ with $\chi \in {\rm Y}(r)$.
Let $\eta : \Gamma \rightarrow \{ \pm 1\}$ be a morphism. \ps
\begin{itemize}
\item[(i)] The pair $(r,\eta)$ satisfies {\rm (U1)} if, and only if,
$(r_{|\Gamma_0},\eta_{|\Gamma_0})$ satisfies {\rm (U1)}.\ps
\item[(ii)] If $(r_{|\Gamma_0},\eta_{|\Gamma_0})$ is unacceptable, 
then $r_{|\Gamma_0}$ is of type {\rm I}.
\end{itemize}
\end{lemma}

\begin{pf} 
The {\it only if} part of (i) is obvious. 
For the {\it if} part, assume $(r_{|\Gamma_0},\eta_{|\Gamma_0})$ satisfies {\rm (U1)}.
We have to show that for $\gamma \in \Gamma \smallsetminus \Gamma_0$, $r(\gamma)$ and $\eta(\gamma)r(\gamma)$ are conjugate in ${\rm Spin}(7)$,
or equivalently (by Proposition~\ref{prop:propequat}), 
that $\eta(\gamma)$ is an eigenvalue of $\gamma$ in $E$.
We conclude as both $1$ and $\chi(\gamma)=-1$ are eigenvalues of $\gamma$ on $E$, 
as $\dim E$ is odd and $\chi \in {\rm Y}(r)$. 
Assertion (ii) is clear since we have $1=\chi_{|\Gamma_0} \in {\rm Y}(r_{|\Gamma_0})$.
\end{pf} 

Recall the set ${\rm E}(r)$ defined in Remark~\ref{rem:etar}.

\begin{coro}
\label{cor:classII} 
Let $r : \Gamma \rightarrow {\rm Spin}(7)$ be a morphism such that $\chi \in {\rm Y}(r)$. 
Then $r$ is unacceptable of type ${\rm II}$ with respect to $\chi$
if, and only if, $r_{|\Gamma_0}$ is unacceptable of type ${\rm I}$ and at least one 
element of ${\rm E}(r_{|\Gamma_0})$ extends to an order $2$ character of $\Gamma$.
\end{coro}

\begin{pf} 
Set $r_0=r_{|\Gamma_0}$.
If $r$ is unacceptable of type II with respect to $\chi$, then by definition 
there is 
$\eta \in {\rm E}(r)$ such that $\eta_{|\Gamma_0} \in {\rm E}(r_0)$ (in particular, $\eta_{|\Gamma_0} \neq 1$). Also, $r_0$ is of type {\rm I} by Corollary~\ref{lem:resGo} (ii).
Conversely, assume $r_0$ is unacceptable of type {\rm I} and 
choose $\eta_0$ in ${\rm E}(r_0)$ extending to $\eta : \Gamma \rightarrow \{\pm 1\}$.
Then $(r,\eta)$ satisfies {\rm (U1)} by Corollary~\ref{lem:resGo} (i), and $(r_0,\eta_0)$
is unacceptable: $r$ is of type {\rm II} with respect to $\chi$.
\end{pf}

This corollary reduces the classification of type {\rm II} unacceptable 
morphisms to that of type I ones, done in the previous section, 
together with some extension problem.
We shall not say more about this problem here, 
and are happy to leave this task to a motivated reader.
We shall content ourselves below with one example of unacceptable type II morphism, 
and with two criteria for their inexistence.

\begin{example} 
\label{ex:typeII}
Let $\vartheta \in N \smallsetminus N^0$ be an order $2$ element with $\vartheta g \vartheta^{-1} = \overline{g}$ for all $g \in {\rm SU}(4)$. For all $(g_1,g_2) \in {\mathcal{G}} \subset {\rm O}(2)^{\pm } \times {\rm O}(2)^{\pm }$ we have $\vartheta (g_1,g_2) \vartheta^{-1} = (\pm g_1,\pm g_2)$. 
Since the $4$ characters of ${\rm O}(2)^{\pm}$ are trivial on $-{\rm 1}_2$, 
the character $\upsilon$ of ${\mathcal{G}}$ extends to ${\mathcal{G}}' = {\mathcal{G}} \rtimes \langle \vartheta \rangle$, and with order $2$ as $\vartheta^2=1$. By Proposition~\ref{prop:ex2Gamma} and Corollary~\ref{cor:classII},
the inclusion of ${\mathcal{G}}'$ in ${\rm Spin}(7)$ is unacceptable of type {\rm II} with respect to 
$\kappa_{|{\mathcal{G}}'}$.
\end{example} 


\begin{prop}
\label{prop:critnotIIa}
If $r : \Gamma \rightarrow {\rm Spin}(7)$ is unacceptable of type {\rm II}, 
then the group $\Gamma$
has a quotient isomorphic to $(\Z/2\Z)^3$ or to $\Z/2\Z \times \Z/4\Z$.
\end{prop}


\begin{pf} 
Choose $\chi$ such that $r$ is of type {\rm II} with respect to $\chi$.
Set $\Gamma_0=\ker\, \chi$ and $r_0 = r_{|\Gamma_0}$.
Choose $\eta \in {\rm E}(r)$ 
such $(r_0,\eta_{|\Gamma_0})$ is unacceptable. 
By Corollary~\ref{cor:uniquedelta}, there is a unique non trivial character $\delta_0$ in ${\rm Y}(r_0)$.
But the $7$-dimensional representation $\pi \circ r_0$ of $\Gamma_0$ extends to $\Gamma$, 
so it is isomorphic to its outer conjugate under
$\Gamma/\Gamma_0 \simeq \Z/2\Z$, hence so is $\delta_0$. It follows that $\delta_0$ extends to a character
$\delta : \Gamma \rightarrow \mu_4$. 
Consider the morphism $f: \Gamma \rightarrow \{ \pm 1\}^2 \times \mu_4$, given by $f=(\eta,\chi,\delta)$. As $\eta$ and $\delta_0$ are distinct of order $2$ over $\Gamma_0$,  
and $\chi \neq 1$, the five characters $1,\chi, \eta, \delta,\eta\delta$ of $\Gamma$ are distinct. We have thus 
$f(\Gamma_0)\simeq (\Z/2\Z)^2$, $|f(\Gamma)|=2|f(\Gamma_0)|=8$,  
$f(\Gamma) \simeq \Z/2\Z \times \Z/4\Z$ if $\delta$ has order $4$,  
and $f(\Gamma) \simeq (\Z/2\Z)^3$ otherwise. 
\end{pf}

\begin{definition}
\label{def:discrete}
A group morphism $r : \Gamma \rightarrow G$ is called {\it discrete} 
if the centralizer of ${\rm Im} \,r$ in $G$ is finite. 
\end{definition}

In the case $G$ is ${\rm Spin}(n)$, it is equivalent 
to ask that $\pi \circ r : \Gamma \rightarrow {\rm SO}(n)$ is discrete,
by arguments given in the proof of Proposition~\ref{prop:critU2}. 
Type I unacceptable morphisms are not discrete by Formula~\eqref{eq:decCGF}. 
On the other hand, the inclusion $\mathcal{G}' \rightarrow {\rm Spin}(7)$ 
in Example~\ref{ex:typeII} is discrete. 
Indeed, its centralizer must be included in $N$, 
hence coincides with $\mu_2 \times \mu_2$.

\begin{prop} 
\label{prop:critnotIIb}
If $r : \Gamma \rightarrow {\rm Spin}(7)$ is discrete and unacceptable of type {\rm II}, 
then the group $\Gamma$ has a quotient isomorphic to $(\Z/2\Z)^3$.
\end{prop}

\begin{pf} 
Define $\chi,\eta,r_0, \delta_0, \delta$ and $f$ as the proof of Proposition~\ref{prop:critnotIIa}.
By Corollary~\ref{cor:uniquedelta}, the isotypic component $U$ of $\delta_0$ in the $\R[\Gamma_0]$-module $E$ has dimension $2$.
This plane $U \subset E$ is $\Gamma$-stable as we have ${\rm Y}(r_0)=\{\delta_0\}$.
As $r$ is discrete, the centralizer of the image of $\Gamma$ in ${\rm O}(U)$ is finite,
hence this image is not included in ${\rm SO}(U)$ (infinite abelian).
If we write $\Gamma = \Gamma_0 \coprod z \Gamma_0$, 
this forces $z$ to act on $U$ as a reflection. 
But then we have $U \simeq \alpha \oplus \beta$ for some characters $\alpha ,\beta : \Gamma \rightarrow \{\pm 1\}$,
and both $\alpha$ and $\beta$ extend $\delta_0$. 
So we have $\delta \in \{\alpha,\kappa \alpha\}$, $\delta(\Gamma) = \mu_2$,
and by the proof of Proposition~\ref{prop:critnotIIa}, $f(\Gamma) \simeq (\Z/2\Z)^3$.
\end{pf}

\section{Type III unacceptable morphisms}
\label{sect:typeIII}

In this section, we still assume that $\Gamma$ is a group, 
that $\chi : \Gamma \rightarrow \{ \pm 1\}$ is a non trivial character, and we 
denote by $\Gamma_0$ the kernel of $\chi$, an index $2$ subgroup of $\Gamma$. 
\ps

\begin{prop} 
\label{prop:typIII}
Let $r : \Gamma \rightarrow {\rm Spin}(7)$ be a morphism with $\chi \in {\rm Y}(r)$.
If $\eta : \Gamma \rightarrow \{ \pm 1\}$ is an order $2$ character, the following are equivalent:\ps
\begin{itemize}
\item[(i)] $(r,\eta)$ is unacceptable and $(r_{|\Gamma_0},\eta_{|\Gamma_0|})$ is acceptable,\ps
\item[(ii)] we have $\eta \not \in {\rm X}(r)$ but $\eta_{|\Gamma_0} \in {\rm X}(r_{|\Gamma_0})$.\ps
\end{itemize}
\end{prop}

\begin{pf} 
By Lemma~\ref{lem:resGo} (i), the assumption $\eta_{|\Gamma_0} \in {\rm X}(r_{|\Gamma_0})$
implies that $(r,\eta)$ satisfies {\rm (U1)}. The equivalence between (i) and (ii) 
follows then from Definition~\ref{def:Xr}.
\end{pf}

Observe that if $(r,\eta)$ is unacceptable of type {\rm I},
and if $\chi$ is the unique non trivial element in ${\rm Y}(r)$ 
(Corollary~\ref{cor:uniquedelta}), then $(r_{|\Gamma_0},\eta_{|\Gamma_0})$
is acceptable (otherwise it would have type {\rm I} and contain $3$ times the trivial character). 
This is why we exclude this case in the following definition.

\begin{definition}
\label{def:typeIII}
A morphism $r : \Gamma \rightarrow {\rm Spin}(7)$ is called 
{\rm unacceptable of type III} with respect to $\chi$ if we have $1 \notin {\rm Y}(r)$
and if there exists an order two character $\eta$ of $\Gamma$ 
satisfying the equivalent properties (i) and (ii) in Proposition~\ref{prop:typIII}.
For such an $\eta$ we also say that $(r,\eta)$ is (unacceptable) of type {\rm III} with respect to $\chi$.
\end{definition}
\ps

By definition, any unacceptable ${\rm Spin}(7)$-valued morphism $r$ is either of type {\rm I}, or 
of type {\rm II} or ${\rm III}$ for some element in ${\rm Y}(r)$.
Let us analyze in more details condition (ii) in Proposition~\ref{prop:typIII}.

\ps

\begin{propdefi} 
\label{prop:devtypIII}
Assume $r$ is unacceptable of type {\rm III} with respect to $\chi$ and
write $E \simeq F \oplus \chi$ as an $\R[\Gamma]$-module. 
There is an order $2$ character $\eta$ of $\Gamma$ with $\eta \notin {\rm X}(r)$,
and an $\R[\Gamma_0]$-submodule $V_0 \subset F$ (not necessarily irreducible)
with determinant $\eta_{|\Gamma_0}$, such that one 
of the following holds: \ps\ps
\begin{itemize}
\item[(a)] either $\dim V_0=2$ and there is an $\R[\Gamma]$-module $Q$, 
with $\dim Q\,=\,2$ and $\det Q=\chi$,  satisfying 
$F \simeq Q \oplus {\rm Ind}_{\Gamma_0}^\Gamma V_0$,\ps
\item[(b)] or $\dim V_0=3$ and we have an $\R[\Gamma]$-module isomorphism
$F \simeq {\rm Ind}_{\Gamma_0}^\Gamma V_0$.\ps
\end{itemize}
If we are in case {\rm (a)} {\rm (}resp. {\rm (b)}{\rm )}, we say that $(r,\eta)$ is unacceptable {\rm of type} 
{\rm IIIa} {\rm (}resp. {\rm IIIb}{\rm )} with respect to $\chi$.
\end{propdefi}

\begin{pf}
By the Definition~\ref{def:typeIII}, there is an order $2$ character
 $\eta'$ of $\Gamma$ and  an $\R[\Gamma_0]$-submodule $U_0 \subset F$ with $\det U_0 = \eta'_{|\Gamma_0}$
 and $\eta' \notin {\rm X}(r)$.
Note that $U_0$ is not $\Gamma$-stable, otherwise the character $\det U_0 \in {\rm X}(r)$ 
would coincide with $\eta'$ on $\Gamma_0$, which forces $\det U_0 \,=\,\eta'$ or $\chi \eta'$, hence 
$\eta' \in {\rm X}(r)$. The biggest $\Gamma$-stable subspace in $U_0$ is 
$U := U_0 \cap zU_0$, where $z$ is any element in  $\Gamma \smallsetminus \Gamma_0$. 
We have thus an $\R[\Gamma_0]$-module decomposition 
$U_0 = U \perp V_0$, with $V_0 \neq 0$, and $\det \,V_0\, = \,(\eta' \det U)_{|\Gamma_0}$.
The character $\eta:=\eta' \det U$ is not in the group ${\rm X}(r)$ (see also Remark~\ref{rem:etar}),
and satisfies $\det\,V_0= \eta_{|\Gamma_0}$.
We also have $V_0 \cap z V_0 =\{0\}$ by construction, and thus an $\R[\Gamma]$-module
embedding ${\rm Ind}_{\Gamma_0}^\Gamma \,V_0\, \hookrightarrow F$, as well as
$1 \leq \dim V_0 \leq 3$. \ps
Assume first $\dim V_0$ is odd and set $W_0 = \eta_{|\Gamma_0} \otimes V_0$.
We have $\det W_0 = 1$ and an isomorphism 
${\rm Ind}_{\Gamma_0}^\Gamma V_0 \simeq 
\eta \otimes {\rm Ind}_{\Gamma_0}^\Gamma W_0$.
In the case $\dim W_0=1$ we must have $W_0 \simeq 1$, hence
${\rm Ind}_{\Gamma_0}^\Gamma W_0 \simeq 1 \oplus \chi$. 
But this shows $\eta \in {\rm Y}(r) \subset {\rm X}(r)$, 
a contradiction. 
This proves $\dim V_0=3$ and we are in case (b).\ps
Assume now $\dim V_0=2$. Define $Q$ as the orthogonal of 
$P={\rm Ind}_{\Gamma_0}^\Gamma \,V_0$
in $F$. 
By Lemma~\ref{lem:trancar} (ii), we have
 $\det P = \chi^2 t=t$ 
where $t : \Gamma \rightarrow \{ \pm 1\}$ is the transfer of 
$\det V_0 = \eta_{|\Gamma_0}$ to $\Gamma$. As we have $\eta^2=1$ 
and $\eta$ is a character of $\Gamma$, we have $t=1$ by Lemma~\ref{lem:trancar} (i). 
This shows $\det P=1$, hence $\det Q = \det F = \chi$.
\end{pf}


\noindent 

\begin{center} {\sc Example 3} \end{center}
\ps\ps
\newcommand{\uH}{\underline{\mathcal{H}}}
Fix a decomposition $F = P \perp P' \perp Q$ with $P,P'$ and $Q$ 
each of dimension $2$. Define $\uH \subset {\rm O}(F)$ as the subgroup 
of isometries that:\ps\ps
-- preserve the pair $\{P,P'\}$ and $Q$,\ps
-- have determinant $1$ on $P \oplus P'$, and\ps
-- whose signature on $\{P,P'\}$
coincides with their determinant on $Q$. \ps\ps
 \noindent Define $\uH_0$ as the subgroup of $\uH$ with trivial determinant on $Q$. 
We clearly have $$\uH_0 \simeq ({\rm O}(2) \times {\rm O}(2))^{\det_1 = \det_2} \times {\rm SO}(2)$$ 
and using this identification $\uH$ may be identified with the semidirect product of $\Z/2\Z$ 
by $\uH_0$ with respect to the involution $(g_1,g_2,g_3) \mapsto (g_2,g_1,g_3^{-1})$. 
Set $\alpha = \det Q$ and denote by $\varepsilon$ any of the two characters of $\uH$ that coincide 
with $\det P=\det P'$ on $\uH_0$ (the other one being $\alpha \varepsilon$). 
Define $\mathcal{H} \subset N \subset {\rm Spin}(7)$ 
as the inverse image of $\uH$.
The characters $\alpha$ and $\varepsilon$ may be viewed as characters of $\mathcal{H}$ by inflation, 
and we set $\mathcal{H}_0 = \ker \alpha$. 
Consider the element $n \in N \smallsetminus N^0$ acting trivially 
on $P'$ and $Q$, and as a reflection on $P$. Observe that 
$n$ normalizes $\mathcal{H}$ by preserving $\mathcal{H}_0$ 
(hence the character $\alpha$) and exchanges
$\varepsilon$ and $\alpha \varepsilon$, as we have $\alpha(\theta n \theta^{-1} n^{-1})=-1$.
We have $F \simeq \,Q \oplus  {\rm Ind}_{\mathcal{H}_0}^\mathcal{H} P$ as $\R[\mathcal{H}]$-modules.
We denote by $\rho : \mathcal{H} \rightarrow {\rm Spin}(7)$ the natural inclusion.

\begin{thm} 
\label{thm:exclasstypIIIa}
The pair $(\rho,\varepsilon)$ is unacceptable of type
{\rm IIIa} with respect to $\alpha$. 
Conversely, if we have $r : \Gamma \rightarrow {\rm Spin}(7)$ 
and $\eta : \Gamma \rightarrow \{ \pm 1\}$
with $(r,\eta)$ unacceptable of type {\rm IIIa} with respect to $\chi$,
then up to replacing $r$ with a conjugate
there is a morphism $f : \Gamma \rightarrow \mathcal{H}$ 
with $\rho \circ f = r$, $\chi = \alpha \circ f$ and $\eta = \varepsilon \circ f$.
\end{thm}

\ps\ps

\begin{pf} 
The characters $\alpha$ and $\varepsilon$ of $\mathcal{H}$ are distinct and non trivial.
Set $\rho_0=\rho_{|\mathcal{H}_0}$ and $\varepsilon_0=\varepsilon_{|\mathcal{H}_0}$.
By Proposition~\ref{prop:critU2}, the morphisms $\rho_0$ and $\varepsilon_0 \rho_0$ are 
${\rm Spin}(7)$-conjugate, since ${\rm X}(\rho_0)$ contains $\varepsilon_0=\det P$.
 That $(\rho,\varepsilon)$ satisfies (U1) follows then from Lemma~\ref{lem:resGo} (i).
The irreducible sub-representations of the $\R[\mathcal{H}]$-module $F$
are clearly $Q$ and ${\rm Ind}_{\mathcal{H}_0}^{\mathcal{H}} P$,
whose determinant are respectively $\alpha$ and $1$.
We have thus ${\rm X}(\rho)=\{1,\alpha\}$ by Proposition~\ref{prop:critU2},  
${\rm Y}(\rho)=\{\alpha\}$,
and $(\rho,\varepsilon)$ 
satisfies (U2) as well since $\varepsilon \notin {\rm X}(\rho)$ : it is unacceptable,
of type IIIa with respect to $\alpha$ by construction.\ps

Suppose conversely that $(r,\eta)$ is unacceptable of type IIIa with respect to $\chi$. 
Choose $z \in \Gamma\smallsetminus \Gamma_0$.
By Lemmas~\ref{lem:propF} (ii) and~\ref{lem:realind},
there is an $\R[\Gamma_0]$-submodule $V_0 \subset F$ 
of dimension $2$ and determinant $\eta_{|\Gamma_0}$
with $F = V_0 \perp zV_0 \perp Q'$, and $Q'$ is $\Gamma$-stable 
with determinant $\chi$. Choose $g \in {\rm Spin}(7)$ with 
$g(V_0)=P$, $g(zV_0)=P'$ and $g(Q')=Q$, we have 
$g\, r(\Gamma) \,g^{-1} \subset \rho(\mathcal{H})$.
Up to replacing $r$ with its $g$ conjugate, 
the morphism $f = \rho^{-1} \circ r$ satisfies
$\rho \circ f = r$, $\chi = \alpha \circ f$ and $\eta_{|\Gamma_0} = \varepsilon \circ f_{|\Gamma_0}$.
The last condition implies $\varepsilon \circ f = \eta$ or $\eta \chi$.
So either $f$, or the conjugate of $f$ by the element $n$ defined above, has all the required properties.
\end{pf}

\begin{remark}
\label{rem:typIIIa}
Of course, if we have a morphism $f : \Gamma \rightarrow \mathcal{H}$, 
and if we set $r= \rho \circ f$, $\chi = \alpha \circ f$ and $\eta = \varepsilon \circ f$,
then $(r,\eta)$ satisfies {\rm (U1)} but not automatically {\rm (U2)}, 
as we may have $\eta \in {\rm X}(r)$ even if
$\varepsilon \notin {\rm X}(\rho)$. We leave to the reader a study of the condition {\rm (U2)} 
in the spirit of that made in Proposition~\ref{prop:classnonex}.
\end{remark}
\ps

\begin{remark}
{\rm ({\it An alternative description of $\mathcal{H}$})
\label{rem:altH}
By definition, we have 
$$\mathcal{H} \subset {\rm Spin}(A)\cdot {\rm Spin}(B), 
\,\,{\rm with}\, A=L \perp Q\,\, {\rm and} \,\,B = P \perp P',$$
and thus ${\rm Spin}(A) \simeq {\rm SU}(2)$ and ${\rm Spin}(B) \simeq {\rm SU}(2) \times {\rm SU}(2)$.
The neutral component $\mathcal{H}^0$ of $\mathcal{H}$ is 
a maximal torus of ${\rm Spin}(A)\cdot {\rm Spin}(B)$. 
Let $T$ be the diagonal torus of ${\rm SU}(2)$, $C$ the normalizer of $T$ in ${\rm SU}(2)$ 
and $s : C \rightarrow \{ \pm 1\}$ the order $2$ character with kernel $T$. There is a morphism 
$\xi : {\rm SU}(2)^3 \rightarrow {\rm Spin}(7)$ with $\xi({\rm SU}(2) \times 1)={\rm Spin}(A)$, 
$\xi(1 \times {\rm SU}(2)^2)={\rm Spin}(B)$ and $\xi(T^3)=\mathcal{H}^0$. 
We have $\ker \xi \,= \,\langle (-{\rm 1}_2,-{\rm 1}_2,-{\rm 1}_2) \rangle \simeq \Z/2\Z$. 
The group $\mathcal{H}$ has index $2$ in $\xi(C^3)$ and a simple computation
shows 
\begin{equation}
\label{eq:calHH} 
\mathcal{H}\,=\,\xi(H) \, \, \,\,{\rm with}\, \, \,\,H=\{ (c_1,c_2,c_3) \in C^3 \, |\, \, s(c_1)s(c_2)s(c_3)=1\}.
\end{equation}
Moreover, if $s_i : H \rightarrow \{\pm 1\}$ denotes the character $(c_1,c_2,c_3) \mapsto s(c_i)$,
we have $\alpha \circ \xi = s_1$ and $\varepsilon \circ \xi = s_2$ or $s_3$. 
}
\end{remark}

\ps\ps
\begin{center} {\sc Example 4} \end{center}
\ps\ps
\newcommand{\uG}{\underline{\mathcal{I}}}

Consider an orthogonal decomposition $F = T \perp T'$ with $\dim T = \dim T' = 3$.
Define $\uG \subset {\rm O}(F)$ as the subgroup 
of isometries that:\ps\ps
-- preserve the pair $\{T,T'\}$, and\ps\ps
-- whose signature on $\{T,T'\}$
coincides with their determinant on $F$. \ps\ps
\noindent Its neutral component is $\uG^0 = {\rm SO}(T) \times {\rm SO}(T') \simeq {\rm SO}(3) \times {\rm SO}(3)$ and we have 
\begin{equation} \label{eq:strucuG} \uG/\uG^0 \simeq \Z/2\Z \times \Z/2\Z.\end{equation}
Indeed, $\uG$ is generated by any order $2$ element in ${\rm O}(F)$ exchanging $T$ and $T'$ (those 
elements have determinant $-1$) and by its index $2$ subgroup 
$$\uG_0=\{ (g_1,g_2) \in {\rm O}(T) \times {\rm O}(T') \,\,|\,\, \det g_1 = \det g_2\} \simeq ({\rm O}(3) \times {\rm O}(3))^{\det_1=\det_2}$$ 
Set $\alpha = \det F$, so that we also have $\uG_0=\ker \alpha$.
Denote by $\varepsilon$ any of the two characters of $\uG$ that coincide 
with $\det T=\det T'$ on $\uG_0$ (the other one being then $\alpha \varepsilon$);
we have $\varepsilon^2 =1$ by \eqref{eq:strucuG}.
Define $\mathcal{I} \subset N \subset {\rm Spin}(7)$ 
as the inverse image of $\uG$.
The characters $\alpha$ and $\varepsilon$ 
may be viewed as characters of $\mathcal{I}$ by inflation, 
and we set $\mathcal{I}_0 = \ker \alpha$. 
We also have $F \simeq \, {\rm Ind}_{\mathcal{I}_0}^\mathcal{I} T$ 
as $\R[\mathcal{I}]$-modules (irreducible). 
Denote by $\rho : \mathcal{I} \rightarrow {\rm Spin}(7)$ the natural inclusion.
The proof of the following theorem is similar to that of Theorem~\ref{thm:exclasstypIIIa}.

\begin{thm} 
\label{thm:exclasstypIIIb}
The pair $(\rho,\varepsilon)$ is unacceptable 
of type {\rm IIIb} with respect to $\alpha$. 
Conversely, if we have $r : \Gamma \rightarrow {\rm Spin}(7)$ 
and $\eta : \Gamma \rightarrow \{ \pm 1\}$
with $(r,\eta)$ unacceptable of type {\rm IIIb} with respect to $\chi$,
then up to replacing $r$ with a conjugate
there is a morphism $f : \Gamma \rightarrow \mathcal{I}$ 
with $r= \rho \circ f$, $\chi = \alpha \circ f$ and $\eta = \varepsilon \circ f$.
\end{thm}

We end this section with several remarks.
First, an analogue of Remark~\ref{rem:typIIIa} also applies in this context.
Moreover:

\begin{remark}
\label{rem:solidplat}
The pair $(\rho,\varepsilon)$ in Theorem~\ref{thm:exclasstypIIIb} is an unacceptable morphism with biggest possible image, of dimension $6$.
Using that theorem and the well-known classification of closed subgroups of ${\rm SO}(3)$, we easily find examples of unacceptable ${\rm Spin}(7)$-valued morphisms $r$ 
of type {\rm IIIb} such that $\pi \circ r(\Gamma)$ is an extension of $\Z/2\Z \times \Z/2\Z$ by $H \times H$
with $H \simeq {\rm A}_4$, ${\rm S}_4$ or ${\rm A}_5$. 
\end{remark}

We leave as an exercise to the reader to verify that the following alternative description
of $\mathcal{I}$ holds.

\begin{remark}
{\rm ({\it An alternative description of $\mathcal{I}$})
There are closed normal subgroups $\mu$ and ${\rm O}$ of $\mathcal{I}$,  
with $\mu \simeq \mu_4$ and ${\rm O} \simeq {\rm O}(4)$, as well as a decomposition
\begin{equation}
\mathcal{I} \,= \,\mu \cdot {\rm O}\, \, {\rm with}\,\, \mu \cap {\rm O}=\{\pm 1\},
\end{equation}
such that if $i$ is a generator of $\mu$,  
then for all $g \in {\rm O}$ we have
$$ g\,i \,g^{-1} \,=\, \det (g)\,i.$$ 
We also have $\alpha_{|\mu}=1$, $\alpha_{|{\rm O}}=\det$ and $\varepsilon_{|\mu} \neq 1$.
}
\end{remark}

\ps

\begin{remark} 
{\rm 
({\it Morphisms of both types {\rm IIIa} and {\rm IIIb}})
We mention that there is an example of a triple $(r,\eta,\chi)$ such that $(r,\eta)$
is unacceptable of both types {\rm IIIa} and {\rm IIIb} with respect to $\chi$. 
Furthermore, this example is unique in a natural sense. We omit the details, 
but simply say that in this example, we have $\pi (r(\Gamma)) \simeq {\rm D}_8 \times \mu_2$
and the $\R[\Gamma]$-module $F$ is isomorphic to 
$P \,\oplus\, \,\eta \otimes P\, \,\oplus \,\,\det P\, \,\oplus \, \,\chi  \det P$,
with $P$ irreducible and $\dim P=2$.}
\end{remark}
 
\section{A few examples and properties in the Weil group case}
\label{sec:weilcase}

In this section, $p$ is a prime, $F$ denotes a finite extension of the field $\Q_p$ of $p$-adic numbers,
and $q$ is the cardinality of the residue field of $F$ (a $p$-th power). 
We denote by ${\rm W}_F$ the Weil group of $F$ (see~\cite{tate}), a locally compact topological group, 
and we consider the unacceptable continuous group morphisms ${\rm W}_F \rightarrow {\rm Spin}(7)$. As already mentioned in the introduction, 
these morphisms are of interest for the approach toward a local Langlands correspondence for the group ${\rm PGSp}_6(F)$ discussed in~\cite{gansavin}. 
We start with a simple lemma.
\ps
\begin{lemma} 
\label{lem:conteta}
Assume $\Gamma$ is a topological group and $r : \Gamma \rightarrow {\rm Spin}(n)$ is continuous.
Then the elements of ${\rm X}(r)$ and of ${\rm E}(r)$ are continuous characters of $\Gamma$.
\end{lemma}
\ps
\begin{pf} 
Recall the morphism $\pi : {\rm Spin}(n) \rightarrow {\rm SO}(n)$.
The assertion about ${\rm X}(r)$ is clear as $\pi$ is continuous. 
The one about $\eta$ follows from the fact that there 
is an open neighborhood $U$ of $1 \in {\rm Spin}(n)$ such that for all $g \in U$, $-g$ is not conjugate to 
$g$. Indeed, there is a neighborhood $U$ of $1$ satisfying ${\rm id}_E+\pi(U) \subset {\rm GL}(E)$ and 
we conclude by Lemma~\ref{lem:critamb} (or Remark~\ref{rem:corprop} for $n$ even).
\end{pf}
\ps

That being said, we only consider {\it continuous} morphisms or characters from now on, without further mention.
As is well-known, we have 
{\scriptsize
\begin{equation}
\label{eq:deltaF}
{\rm Hom}({\rm W}_F,\{\pm 1\}) \simeq F^\times/F^{\times,2} \simeq (\Z/2\Z)^{\delta_F} \, \, \, {\rm with}\, \,\,
\delta_F = \left \{ \begin{array}{cl} 
$2$ & {\rm if} \,\,p\,\,{\rm is\,\, odd}, \\
2+[F:\Q_2] & {\rm for}\,\, p=2. 
\end{array}\right.
\end{equation}
}
In particular, we always have $\delta_F\geq 2$, a necessary condition for the existence of unacceptable
morphisms by Corollary~\ref{coro:homz22}.\footnote{For the Archimedean local fields $F=\R$ or $F=\C$,
this same corollary also shows that any continuous morphism ${\rm W}_F \rightarrow {\rm Spin}(7)$ is acceptable.
} If $E$ is a finite extension of $F$, we will denote by 
${\rm N}_{E/F} : E^\times \rightarrow F^\times$ the norm morphism and by 
${\rm S}^1(E/F)$ its kernel, a compact subgroup of $E^\times$.\ps

\begin{lemma} 
\label{lem:prolcar}
Let $E$ be a quadratic extension of $F$ and $\epsilon : F^\times \rightarrow \C^\times$ a character. 
For any integer $n\geq 1$,
there is a character $E^\times \rightarrow \C^\times$ whose restriction to $F^\times$ is $\epsilon$ and whose 
restriction to ${\rm S}^1(E/F) \subset E^\times$ has order divisible by $p^n$.
\end{lemma}

\begin{pf} 
Set $S={\rm S}^1(E/F)$.
The subgroup $F^\times \cdot S$ of $E^\times$ is open, of finite index, 
and we have $F^\times \cap S = \{ \pm 1\}$. 
We may write $S = S_{\rm tor} \cdot S_f$ with $S_{\rm tor} \subset S$
the finite torsion subgroup and $S_f \simeq \Z_p^d$ with $d=[F:\Q_p]$.
We extend first $\epsilon$ to a character $\epsilon'$ of $F^\times \cdot S_{\rm tor}$.
As we have $(F^\times \cdot S_{\rm tor}) \cap S_f =\{1\}$, we may extend $\epsilon'$ to a character
$\epsilon''$ of $F^\times \cdot S$ so that $\epsilon''_{|S_f}$ has order $p^n$.
Any extension of $\epsilon''$ to $E^\times$ does the trick.
\end{pf}

Recall the {\it reciprocity isomorphism} ${\rm rec}_F: F^\times \isomo {\rm W}_F^{\rm ab}$ from local 
class field theory \cite[\S 1.1]{tate}. For any quadratic extension $E/F$ we denote by 
${\rm sgn}_{E/F}: {\rm W}_F \rightarrow \{ \pm 1\}$ the order $2$ character with kernel ${\rm W}_E$. 
The kernel of ${\rm sgn}_{E/F} \circ {\rm rec}_F$ is the index $2$ subgroup 
${\rm N}_{E/F}(E^\times)$ of $F^\times$. 
If $c : E^\times \rightarrow \C^\times$ is a character, we consider the induced representation 
${\rm I}(c) \,:=\, {\rm Ind}_{{\rm W}_E}^{{\rm W}_F}\, \, c \circ {\rm rec}_E^{-1}$.
For general reasons we have 
\begin{equation}
\label{eq:inducedic}
\, \, {\rm I}(c) \simeq {\rm I}(c)^\ast \otimes \det {\rm I}(c)\,\,\,\,{\rm and}\, \, \,\,{\rm I}(c) \simeq {\rm I}(c) \otimes {\rm sgn}_{E/F}.
\end{equation}
By~\cite[\S 1.1 (W3)]{tate} the transfer of $c \circ {\rm rec}_E^{-1}$ 
to ${\rm W}_F$ is $c_{|F^\times} \circ {\rm rec}_F^{-1}$, so we have
\begin{equation}
\label{eq:detic}
\det {\rm I}(c)\, \circ {\rm rec}_F\, =\, {\rm sgn}_{E/F} \circ {\rm rec}_F \,\cdot \,c_{|F^\times}
\end{equation}
by Lemma~\ref{lem:trancar} (ii).
\ps\ps
\begin{prop} 
\label{prop:infinitetypeI}
There exist continuous morphisms $r : {\rm W}_F \rightarrow {\rm Spin}(7)$ which are unacceptable of type {\rm I} and of arbitrary large finite image.
\end{prop}
\ps
\begin{pf} 
By~\eqref{eq:deltaF} we may choose two different quadratic extensions $E$ and $E'$ of $F$. 
Set $s= {\rm sgn}_{E/F}$ and $s'={\rm sgn}_{E'/F}$, two order $2$ characters of ${\rm W}_F$.
The character $s'':=ss'$ is the character ${\rm sgn}_{E''/F}$ of the third quadratic extension $E''$ of $F$ in
the compositum $E \cdot E'$.
By Lemma~\ref{lem:prolcar}, we may choose a character $c : E^\times \rightarrow \C^\times$ 
with $c_{|F^\times} =s' \circ {\rm rec}_F$ and of arbitrary large finite order over ${\rm S}^1(E/F)$.
Such a character has finite image as $s'$ has this property and $F^\times \cdot {\rm S}^1(E/F)$ has finite
index in $E^\times$.\par
By~\eqref{eq:detic} we have $\det {\rm I}(c) = ss' =s''$.
By~\eqref{eq:inducedic} we deduce ${\rm I}(c)^\ast \simeq {\rm I}(c) \otimes s'$. 
Using $s' \neq \det {\rm I}(c)$, Lemma~\ref{lem:D8a} (i) shows that ${\rm I}(c)$ defines a morphism $r : {\rm W}_F \rightarrow {\rm O}(2)^{\pm }$
with $\mu \circ r = s'$, $\det \circ r = s''$ and thus $\epsilon \circ r = s$.
The group $r({\rm W}_F)$ is not isomorphic to ${\rm D}_8$ as long as we choose $c$ to have order $>2$.
By exchanging the roles of $E$ and $E'$ we may also find 
a morphism $r' : {\rm W}_F \rightarrow {\rm O}(2)^{\pm }$
with $\mu \circ r = s$, $\det \circ r = s''$, $\epsilon \circ r = s'$ and $r'({\rm W}_F) \not \simeq {\rm D}_8$.
The pair $\rho:=(r,r')$ so defined is a morphism $\rho : {\rm W}_F \rightarrow \mathcal{G}$ 
with $\nu \circ \rho = s$ and ${\rm d} \circ \rho = s''$. By Proposition~\ref{prop:converseN}, the pair $(\rho,s)$ is unacceptable of type I.
As the order of $c$ is arbitrary large, so is the cardinality of the finite group $\rho({\rm W}_F)$.
\end{pf}

	As a consequence, there are always plenty of unacceptable morphisms ${\rm W}_F \rightarrow {\rm Spin}(7)$. An interesting question, from the point of view of the Langlands classification, is the existence of {\it discrete} unacceptable morphisms ${\rm W}_F \rightarrow {\rm Spin}(7)$ (see Definition~\ref{def:discrete}).  
	

\begin{prop} \begin{itemize}
\item[(a)] For any group $\Gamma$, there is no type {\rm I} discrete unacceptable morphism $\Gamma \rightarrow {\rm Spin}(7)$.\ps
\item[(b)] When $p$ is odd, there is no type {\rm II} discrete unacceptable morphism ${\rm W}_F \rightarrow {\rm Spin}(7)$
\end{itemize}
\end{prop}

\begin{pf}
(a) If $r: \Gamma \rightarrow {\rm Spin}(7)$ unacceptable of type I, 
the $\R[\Gamma]$-module $E$ contains some order $2$ character with multiplicity $2$ by~Corollary \ref{cor:uniquedelta}, 
so the centralizer of $r(\Gamma)$ in ${\rm Spin}(7)$ contains some ${\rm Spin}(2) \simeq {\rm S}^1$ (infinite).   \par
(b)
For $p$ odd, there is no surjective morphism $F^\times \rightarrow (\Z/2\Z)^3$ by~\eqref{eq:deltaF}, 
so the second assertion follows from Proposition~\ref{prop:critnotIIb}.
\end{pf}

On the other hand, there are always plenty of unacceptable discrete morphisms ${\rm W}_F \rightarrow {\rm Spin}(7)$
of type III.

\begin{prop} 
\label{prop:infinitetypeI}
There exist discrete continuous morphisms $r : {\rm W}_F \rightarrow {\rm Spin}(7)$ which are unacceptable of type {\rm IIIa} and of arbitrary large finite image.
\end{prop}

\begin{pf} 
Let $K$ be a Galois extension of $F$ with ${\rm Gal}(K/F) \simeq (\Z/2\Z)^2$, and let 
$F_1, F_2$ and $F_3$ be the three quadratic extensions of $F$ inside $K$.
For each $1 \leq i \leq 3$ we may choose by Lemma~\ref{lem:prolcar} a character 
$c_i : F_i^\times \rightarrow \C^\times$ with 
$(c_i)_{|F^\times}={\rm sgn}_{F_i/F} \circ {\rm rec}_F$ and whose restriction to ${\rm S}^1(F_i/F)$
has an arbitrarily high order $>4$. 
As $F^\times\cdot {\rm S}^1(F_i/F)$ has finite index in $F_i^\times$, the image of $c_i$ is finite.
The representation ${\rm I}(c_i)$ of ${\rm W}_F$ has determinant $1$
by Formula~\eqref{eq:detic}, so each $r_i$ defines a morphism
$r_i : {\rm W}_F \rightarrow {\rm SU}(2)$.\par
We now use the description of the group $\mathcal{H}$ given in Remark~\ref{rem:altH}.
We borrow the notations $T, C, H, s, s_i, \xi$ of this remark.
Up to conjugating $r_i$ we may assume $r_i({\rm W}_{F_i}) \subset T$, 
which forces $r_i({\rm W}_F) \subset C$ as we have $c_i^2 \neq 1$. 
The composition ${\rm W}_F \overset{r_i}{\rightarrow} C \overset{s}{\rightarrow} \{ \pm 1\}$
is ${\rm sgn}_{F_i/F}$ by construction, and we have the identity
$${\rm sgn}_{F_1/F}\,\,{\rm sgn}_{F_2/F}\,\,{\rm sgn}_{F_3/F}\,\,=\,\,1,$$
so that $r:=(r_1,r_2,r_3)$ defines a morphism 
$r: {\rm W}_F \rightarrow H$. 
We now consider $\xi \circ r : {\rm W}_F \rightarrow \mathcal{H}$.
We have $\varepsilon \circ \xi \circ r =  {\rm sgn}_{F_2/F}$ and $\alpha  \circ \xi \circ r = {\rm sgn}_{F_1/F}$.
Then $(\xi \circ r, {\rm sgn}_{F_2/F})$ satisfies (U1) by the first assertion of Theorem~\ref{thm:exclasstypIIIa}, 
and if we can show it satisfies (U2), then it will be of type IIIa with respect to ${\rm sgn}_{F_1/F}$.
So it only remains to show that (U2) holds, {\it i.e.} ${\rm sgn}_{F_2/F} \notin {\rm X}(\xi \circ r)$.\par
	The $\R[H]$-module $E$ has an (absolutely) irreducible decomposition of the form $E = {\rm s}_1 \oplus Q \oplus S$ with $\dim Q=2$, $\det Q = s_1$, $\dim S = 4$ and $\det S = 1$ (we use the letter $S$ here instead of $F$ since the later already denotes the local field). As we have ${\rm sgn}_{F_1/F} \neq {\rm sgn}_{F_2/F}$
it is enough to show that $S$ and $Q$ are absolutely irreducible as $\R[{\rm W}_F]$-modules. 
This is clear for $Q$ as we have $\det Q \neq 1$ and the image of ${\rm W}_F$ in ${\rm O}(Q)$ has order $>4$
by assumption on $c_1$. So we now deal with $S$. For each $1\leq i \leq 3$ define a character 
$a_i: K^\times \rightarrow \C^\times$ by $a_i = c_i \circ {\rm N}_{K/F_i}$. By the basic 
properties of the reciprocity morphisms, we have ${\rm I}(c_i)_{|{\rm W}_K} \simeq a_i \oplus a_i^{-1}$,
and then a $\C[{\rm W}_K]$-module isomorphism
$$S \otimes \C \simeq \,a_2 a_3\, \oplus \, (a_2a_3)^{-1}\,\oplus a_2^{-1}a_3\, \oplus a_2 a_3^{-1}.$$
As each of $s_1,s_2$ and $s_3$ is nontrivial over ${\rm W}_F$, the four characters $a_2^{\pm 1}a_3^{\pm 1}$ are conjugate 
under ${\rm W}_F$. It is thus enough to show that they are distinct, or equivalently, that $a_2^2$, $a_3^2$ and $(a_2a_3)^2$ are all nontrivial over $K^\times$. Denote by $\sigma_i \in {\rm Gal}(K/F)$ the order $2$
element fixing pointwise $F_i$. For $x \in K^\times$ we have $a_i(x) = c_i(x \sigma_i(x))$ by definition,
hence ${a_i}_{|F_i^\times}=c_i^2$, which has order $>2$ by assumption on $c_i$.
Last but not least, for $x \in {\rm S}^1(F_2/F) \subset F_2^\times \subset K^\times$ 
we have just seen $a_2(x)=c_2(x)^2$, but we also have $\tau_3(x)=x^{-1}$, hence $a_3(x)=c_3(xx^{-1})=1$, 
and thus the identity $$(a_2a_3)_{|{\rm S}^1(F_2/F)}=(c_2^2)_{|{\rm S}^1(F_2/F)}.$$
We conclude as the character on the right has order $>2$ by assumption.
\end{pf}
\ps
From the point of view of Langlands' theory of endoscopy, a more specific class of 
discrete morphisms $r : {\rm W}_F \rightarrow {\rm Spin}(n)$ is of interest, namely those such
that the centralizer of $r({\rm W}_F)$ in ${\rm Spin}(n)$ is the center of ${\rm Spin}(n)$, hence {\it as small as possible}.
We call {\it stable} such a morphism, a meaningful terminology from the point of view of endoscopy.
It is easy to see that the discrete type III examples of Proposition~\ref{prop:infinitetypeI} are not stable 
(use Lemma~\ref{lem:stabgen} (i) below). 
For $p$ odd, even more is happily true:
\ps
\begin{prop} 
\label{prop:stabacc}
For $p$ odd, there is no discrete, stable and unacceptable morphism $r : {\rm W}_F \rightarrow {\rm Spin}(7)$.\end{prop}
\ps
\begin{pf} 
Set $\Gamma={\rm W}_F$ and let $r : \Gamma \rightarrow {\rm Spin}(7)$ be discrete, stable and unacceptable. By \eqref{eq:deltaF} and Lemma~\ref{lem:stabgen} (ii) below, the $\R[\Gamma]$-module $E$ has at most $2$ irreducible summands.
An immediate inspection of $E$ in types I, II and III (Formula~\eqref{eq:decCGF}, Corollary~\eqref{cor:uniquedelta} and Proposition-Definition~\eqref{prop:devtypIII}) 
shows that the only possibility is that $r$ is of type IIIb. 
More precisely, there are $2$ characters $\chi,\eta : \Gamma \rightarrow \{ \pm 1\}$
as well as real $3$-dimensional, necessarily irreducible, 
representation $V_0$ of the kernel $\Gamma_0$ of $\chi$,
such that we have an $\R[\Gamma]$-module decomposition 
$$E \simeq \chi \oplus {\rm Ind}_{\Gamma_0}^\Gamma V_0.$$
So $V_0$ gives rise to an irreducible representation $\Gamma_0 \rightarrow {\rm O}(3)$,
and after twisting it by $\det V_0$, to an irreducible representation $\Gamma_0 \rightarrow {\rm SO}(3)$.
But we know since Klein that the finite irreducible subgroups of ${\rm SO}(3)$ 
are isomorphic to ${\rm A}_4$, ${\rm S}_4$ and ${\rm A}_5$. 
But it is well-known that none of these groups can be the Galois group of a finite Galois extension of $p$-adic fields with $p$ odd \cite[\S 13]{weil}. 
\end{pf}

\begin{lemma} 
\label{lem:stabgen}
Assume $r : \Gamma \rightarrow {\rm Spin}(n)$ is discrete and stable with $n$ odd.\ps
\begin{itemize}
\item[(i)] There is no $\R[\Gamma]$-submodule $\{0\} \subsetneq V \subsetneq E$ with $\det V = 1$. \ps
\item[(ii)] If $r$ is furthermore unacceptable, then either $E$ has $\leq 2$ irreducible summands, or there is a surjective group morphism $\Gamma \rightarrow (\Z/2\Z)^3$.
\end{itemize}
\end{lemma}

\begin{pf} 
To prove assertion (i), assume we have a nontrivial $\R[\Gamma$]-stable decomposition $E = A \perp B$.
Then we have $\det A=\det B$, and if this character is trivial the group 
$r(\Gamma)$ falls inside the subgroup ${\rm Spin}(A) \cdot {\rm Spin}(B)$ of ${\rm Spin}(A,B)$, 
whose center strictly contains $\{ \pm 1\}$. For assertion (ii), assume we have a $\Gamma$-stable decomposition $E = E_1 \oplus E_2 \oplus E_3$ with $E_i$ nonzero for each $i$.
By assertion (i), the three characters $\det E_i$ are nontrivial, distinct, and of course in ${\rm X}(r)$. 
But by the unacceptability of $r$ there is another order $2$ character $\eta$ of $\Gamma$, with $\eta \notin {\rm X}(r)$. 
\end{pf}

\ps
The reader aware of Langlands' parameterizations knows that we have to consider more generally 
continuous morphisms ${\rm W}_F \times {\rm SU}(2) \rightarrow {\rm Spin}(7)$.
We did not emphasize this extra ${\rm SU}(2)$ earlier because of the following proposition. 
\ps
\begin{prop}
\label{prop:dealSU2} 
Let $W$ be any group and assume $r : W \times {\rm SU}(2) \rightarrow {\rm Spin}(7)$ is a  
morphism whose restriction to ${\rm SU}(2)$ is continuous and nontrivial. Then $r$ is acceptable.
\end{prop}

\ps

\begin{pf} 
Set $\Gamma = W \times {\rm SU}(2)$ and assume $r$ is unacceptable by contradiction.
There are no non trivial continuous morphisms from ${\rm SU}(2)$ to $\mathcal{G}$ or $\mathcal{H}$,
since those two groups have an abelian neutral component.
By Theorems~\ref{thm:univGG} and~\ref{thm:exclasstypIIIa}, 
the morphism $r$ is neither of type I, nor of type IIIa.
So there are characters $\chi,\eta : \Gamma \rightarrow \{ \pm 1\}$
such that $(r,\eta)$ is of type II or IIIb with respect to $\chi$.
A character of $\Gamma$ is necessarily trivial over $1 \times {\rm SU}(2)$. 
So the kernel $\Gamma_0$ of $\chi$ has the form 
$$\Gamma_0= W_0 \times {\rm SU}(2)$$ for some index $2$ subgroup $W_0$ of $W$.  
By Lemma~\ref{lem:resGo} (ii), if $r_{|\Gamma_0}$ is unacceptable then it is of type I,
which contradicts the previous paragraph applied to $r_{|W_0 \times {\rm SU}(2)}$, 
so $(r,\eta)$ is of type IIIb. 
By Definition-Proposition~\ref{prop:devtypIII}, we have $$E \simeq \chi \oplus {\rm Ind}_{\Gamma_0}^\Gamma V_0$$ for some $\R[\Gamma_0]$-module $V_0$ with $\dim V_0=3$ and $\det V_0 =\eta_{|\Gamma_0}$. 
If $1 \times {\rm SU}(2)$ acts trivially on $V_0$, then it acts trivially as well on $E$, a contradiction.
But then it must act absolutely irreducibly on $V_0$, as we have $\dim V_0=3$.
As $W_0 \times 1$ centralizes $1 \times {\rm SU}(2)$, it must act by multiplication by a real character on $V_0$, 
necessarily equal to $\eta_{|\Gamma_0}$ by taking the determinant. So we have $V_0 \simeq \eta_{|W_0} \boxtimes U$, 
and then 
$${\rm Ind}_{\Gamma_0}^\Gamma V_0  \,\simeq\, ({\rm Ind}_{W_0}^W\, \eta_{|W_0}) \boxtimes U \,\simeq\, \eta_{|W} \boxtimes U \oplus (\chi \eta)_{|W} \boxtimes U.$$
The factor $\eta_{|W} \boxtimes U$ has determinant $\eta$, contradicting the unacceptability of $r$.
\end{pf}

\section{The ${\rm GSpin}(n)$ variant}
\label{sec:gspinn}

In this section, $n$ is any integer $\geq 1$ and we explain how the previous results 
can be applied to study the unacceptable ${\rm GSpin}(n)$-valued morphisms. 
Recall that ${\rm GSpin}(n)$ is the compact subroup of ${\rm Cl}(E)^\times$
generated by its unit scalar subgroup $Z\simeq {\rm U}(1)$ and ${\rm Spin}(n)$ (see Sect.~\ref{sec:defspin}). So $Z$ is central in ${\rm GSpin}(n)$ and we have
\begin{equation}
\label{eq:devGSpin}
{\rm GSpin}(n) \,= \,Z\,\cdot\, {\rm Spin}(n)\,\,\,{\rm and}\, \, \, Z \cap {\rm Spin}(n) = \{ \pm 1\}.
\end{equation}
The morphism $\pi : {\rm Spin}(n) \rightarrow {\rm SO}(n)$ defined {\it loc. cit.} extends to a morphism 
${\rm GSpin}(n) \rightarrow {\rm SO}(n)$ with kernel $Z$ and still denoted by $\pi$. 
We fix a morphism $r : \Gamma \rightarrow {\rm GSpin}(n)$ and set
$$\Gamma(r)= \{ (\gamma,\sigma) \in \Gamma \times {\rm Spin}(n) \,\, |\, \, 
\pi(r(\gamma))=\pi(\sigma) \}.$$
The first projection  
$\Gamma(r) \rightarrow \Gamma$ 
is surjective with kernel $\{ (1,\pm 1) \}$;  any morphism $f$ with source $\Gamma$ 
can be thus inflated to a morphism $\widetilde{f}$ with source $\Gamma(r)$ using this surjection.
 We also have two natural morphisms $r_S:  \Gamma(r) \rightarrow {\rm Spin}(n)$  and $r_Z :  \Gamma(r) \rightarrow Z$, defined by
 $r_S(\gamma,\sigma)= \sigma$ and $r_Z(\gamma,\sigma)=r(\gamma) \sigma^{-1}$, and satisfying 
\begin{equation}
\label{eq:rZrS}
\widetilde{r}(g) \,=\, r_Z(g)\,r_S(g),\, \, \, \forall g \in \Gamma(r).
\end{equation}
We denote by ${\rm G}(r)$ the set of morphisms $\Gamma \rightarrow {\rm GSpin}(n)$
which are element conjugate to $r$, and by ${\rm S}(r)$ the set of morphisms $\Gamma(r) \rightarrow {\rm Spin}(n)$ which are element conjugate to $r_S$. Of course, we have $r \in {\rm G}(r)$ and $r_S \in {\rm S}(r)$.
The group ${\rm GSpin}(n)$ naturally acts on ${\rm G}(r)$ and ${\rm S}(r)$ by conjugation.

\begin{prop} 
\label{prop:gspinn}
Let $r : \Gamma \rightarrow {\rm GSpin}(n)$ be a fixed morphism.
There is a natural, ${\rm GSpin}(n)$-equivariant, bijection $b : {\rm G}(r) \isomo {\rm S}(r)$ satisfying $b(r)=r_S$.
In particular, $r$ is unacceptable if, and only if, $r_S$ is unacceptable. 
\end{prop}

Since we have $\widetilde{\pi \circ r} = \pi \circ r_S$ by Formula~\eqref{eq:rZrS}, 
most of what we have done for ${\rm Spin}(n)$-valued morphisms will apply to ${\rm GSpin}(n)$-valued ones by this proposition. We refer to the proof of Proposition~\ref{prop:gspinkmor} for a concrete example.

\begin{pf} We follow a construction in the proof of \cite[Prop.1.4]{larsen1}.
For $r' \in {\rm G}(r)$ and $g=(\gamma,\sigma) \in \Gamma(r)$, we set 
$b(r')(g):=r_Z(g)^{-1} r'(\gamma) \in {\rm GSpin}(n)$. 
This element is ${\rm Spin}(n)$-conjugate to $r_S(g)=r_Z(g)^{-1} r(\gamma)$
by Formulas~\eqref{eq:rZrS} and~\eqref{eq:devGSpin}. 
As $b(r')$ is a group morphism $\Gamma(r) \rightarrow {\rm GSpin}(r)$, we have $b(r') \in {\rm S}(r)$ and the identity
\begin{equation}
\label{eq:rpZrpS}
\widetilde{r'}(g) \,=\, r_Z(g)\,b(r')(g),\, \, \, \forall g \in \Gamma(r).
\end{equation}
Conversely, for any $f \in {\rm S}(r)$ and $g \in \Gamma(r)$ the element 
$r_Z(g) f(g) \in {\rm GSpin}(n)$ is conjugate to $\widetilde{r}(g)=r_Z(g) r_S(g)$ by assumption.
In particular, it is trivial on $(1,-1)$, and so there is a unique $r' \in {\rm G}(r)$ with $b(r')=f$.
We clearly have $b( g r' g^{-1}) = g b(r') g^{-1}$ for all $r' \in {\rm G}(r)$ and $g \in {\rm GSpin}(n)$.
The last assertion follows since
 ${\rm GSpin}(n)$-conjugacy and ${\rm Spin}(n)$-conjugacy coincide in ${\rm Spin}(n)$ by 
Formula~\eqref{eq:devGSpin}.
\end{pf}

\section{An application in the non compact case}
\label{sec:noncompact}

In this last section, we fix an {\it algebraically closed field $k$ that embeds into $\C$}.
Our aim is to give some statement about morphisms\footnote{In the applications to Galois representations, such as those in~\cite{kretshin}, ${\rm GSpin}$-valued morphisms are much more common than ${\rm Spin}$-valued morphisms.} 
$\Gamma \rightarrow {\rm GSpin}_n(k)$ in the spirit of~\cite[\S 4 \& 5]{kretshin} and 
that follows from our results.
This forces us to discuss an analogue of the acceptability condition for morphisms to linear algebraic $k$-groups that is useful in practice but which slightly differs from the case of morphisms to compact groups. In the proof of Proposition~\ref{prop:gspinkmor} below, and in Proposition~\ref{prop:decp}, we will see how to pass from a setting to the other.\ps 

Let $G$ be a linear algebraic $k$-group, $\Gamma$ an arbitrary group and $r: \Gamma \rightarrow G$
a group morphism.\footnote{We identify a linear algebraic $k$-group with its group of $k$-points, 
as in~\cite{humphreys} for instance.} We denote by ${\rm Zar}(r)$ the Zariski closure of $r(\Gamma)$ in $G$.
Recall that $r$ is called {\it semi-simple} if ${\rm Zar}(r)$ is reductive.
\footnote{We do not assume that a reductive group is connected.} 
Two elements $g,g'$ of $G$ are called {\it ss-conjugate} if the
semi-simple parts in their Jordan decompositions are conjugate,
or equivalently, if $g$ and $g'$ have the same trace in any algebraic $k$-linear 
representation of $G$.
This equivalence relation on $G$ is thus Zariski-closed in $G \times G$. 
Two morphisms $r, r': \Gamma \rightarrow G$ are called {\it element ss-conjugate} 
if, for each $\gamma \in \Gamma$, $r(\gamma)$ and $r'(\gamma)$ are {\it ss-conjugate}.
Finally, we say that a morphism $r : \Gamma \rightarrow G$ is {\it ss-acceptable} if, for
each $r' : \Gamma \rightarrow G$ element ss-conjugate to $r$, 
then $r'$ is actually $G$-conjugate to $r$. 
The following proposition slightly strengthens~\cite[Prop. 1.7]{larsen1}.

\begin{prop} 
\label{prop:decp}
Let $G_1$ and $G_2$ be complex reductive linear algebraic groups, and let $K_1$ and $K_2$
be maximal compact subgroups of $G_1$ and $G_2$ respectively. Assume 
$r,r' : G_1 \rightarrow G_2$ are two algebraic morphisms
with\footnote{We can always achieve this property by conjugating $r$ and $r'$ in $G_2$.} $r(K_1) \subset K_2$ and
$r'(K_1) \subset K_2$, and consider the two morphisms $r_{|K_1}, r'_{|K_1}$ : $K_1 \rightarrow K_2$. Then $r_{|K_1}$ and $r'_{|K_1}$ are element conjugate (resp. $K_2$-conjugate)
if, and only if, $r$ and $r'$ are element ss-conjugate (resp. $G_2$-conjugate). 
\end{prop}

\begin{pf} 
Consider the Cartan (or {\it polar}) decomposition $G_2=K_2 P$ of the linear reductive Lie group $G_2$ with respect to its maximal compact subgroup $K_2$. Recall that $P$ is a subset of $G_2$ stable by $K_2$-conjugacy such that the multiplication $K_2 \times P \rightarrow G_2$ is bijective.
Assume $g x g^{-1} = y$ with $x,y \in K_2$ and $g \in G_2$, and write $g=pk$ with $k \in K_2$ and $p \in P$. The uniqueness of Cartan decomposition and $k P k^{-1} = P$ show $k x k^{-1} = y$ (and $py=yp$).
As any element of $K_2$ is semi-simple, this shows that $r_{|K_1}$ and $r'_{|K_1}$ are element conjugate (resp. $K_2$-conjugate) if $r$ and $r'$ are element ss-conjugate (resp. $G_2$-conjugate).
The converse follows from the Zariski density of $K_1$ in $G_1$, as
``element ss-conjugacy'' is a Zariski closed relation in $G_2 \times G_2$.
\end{pf}

Consider the quadratic space $E=k^n$ with quadratic form $x_1^2+\cdots+x_n^2$.
As in \S\ref{sec:defspin}, we have associated reductive linear algebraic $k$-groups ${\rm SO}(E)$ 
and ${\rm Spin}_n(k) \subset {\rm GSpin}_n(k) \subset {\rm Cl}(E)^\times$, 
and a natural surjective morphism $\pi : {\rm GSpin}_n(k) \rightarrow {\rm SO}(E)$. 

\begin{prop} 
\label{prop:gspinkmor}
Let $r: \Gamma \rightarrow {\rm GSpin}_n(k)$ be a semi-simple morphism with $n\leq 7$.
In the case $n=7$, we assume that one of the two following conditions hold: 
\begin{itemize}
\item[(i)] The $k[\Gamma]$-module $E$ does not contain any
character $c : \Gamma \rightarrow k^\times$ with $c^2=1$.\ps
\item[(ii)] The multiplicity of the weight $0$ in the $k[\Gamma]$-module $E$ is $\leq 2$ and 
$r(\Gamma)$ contains a non trivial unipotent element.
\end{itemize}
Then $r$ is ss-acceptable.
\end{prop}

The first part of condition (ii) means that for some (hence any) maximal torus $T$ of ${\rm Zar}(r)$, 
the invariants of $T$ in $E$ have dimension $\leq 2$.

\begin{pf} Set $G={\rm GSpin}_n(k)$.
Assume $r' : \Gamma \rightarrow G$ is element ss-conjugate to $r$.
In order to show that $r'$ is $G$-conjugate to $r$, we may assume $k=\C$ by the 
Nullstellensatz, since $k$ embeds in $\C$ by assumption. 
In the style of Remark~\ref{rem:topo}, up to replacing $\Gamma$ by 
the Zariski closure of its image in $r \times r' : \Gamma \rightarrow G \times G$, 
and $r$ and $r'$ by the two projections,
we may assume that $\Gamma$ is a complex reductive linear algebraic group, 
and that $r$ and $r'$ are injective algebraic morphisms.  \par
The group ${\rm GSpin}(n)$ is maximal compact in $G$.
Choose $K$ a maximal compact subgroup of $\Gamma$.
Up to replacing $r$ and $r'$ by some $G$-conjugate if necessary, 
we may assume we have $r(K), r'(K) \subset {\rm GSpin}(n)$,
and consider $r_{|K}$ and $r'_{|K} : K \rightarrow {\rm GSpin}(n)$ 
as in Proposition~\ref{prop:decp}.
By this proposition, $r_{|K}$ and $r'_{|K}$ are element conjugate,
and they are ${\rm GSpin}(n)$-conjugate if, and only if, $r$ and $r'$ are $G$-conjugate.
The acceptability of ${\rm Spin}(n)$ for $n\leq 6$, hence that of ${\rm GSpin}(n)$
by Proposition~\ref{prop:gspinn}, concludes the proof for $n \leq 6$. \par
We may thus assume $n=7$ and that the morphism $f:=r_{|K}$ is unacceptable. 
We apply to this $f : K \rightarrow {\rm GSpin}(7)$ the considerations of \S \ref{sec:gspinn}
and use the notations $K(f)$ and $f_S$ {\it loc. cit.}, with $\Gamma=K$ and $r=f$.
By Proposition~\ref{prop:gspinn}, the morphism $f_S : K(f) \rightarrow {\rm Spin}(7)$ 
is unacceptable as well, and $\pi \circ f_S$ factors through $K(f) \rightarrow K$ and 
coincides then with $\pi \circ r_{|K}$. By Theorem~\ref{thm:exceptionPin6},
there is a line $L$ in the Euclidean space $E_\R = \R^7$ 
(a real structure of the complex quadratic space $E$)
on which $K$ acts by a character $K \rightarrow \{ \pm 1\}$.
As $K$ is Zariski dense in $\Gamma$, it follows that $\Gamma$ also acts on $L \otimes \C \subset E$
by such a character, contradicting assumption (i).
So we may assume that (ii) holds.\par
Fix $T$ a maximal torus in $K$.
If $f_S$ is of type I, II or IIIb, it follows from 
Theorem~\ref{thm:univGG}, Corollary~\ref{cor:classII} 
and Theorem~\ref{thm:exclasstypIIIb} respectively 
that the invariants of $T$ in $E_\R$ have dimension $\geq 3$
(consider the action of a maximal torus of $\mathcal{G}$ and $\mathcal{I}$ on $E_\R$).
As $T$ is Zariski dense in a maximal torus of $\Gamma$, 
and as $r$ induces an isomorphism $\Gamma \simeq {\rm Z}(r)$, this contradicts 
the first assertion in assumption (ii).\par
So we may assume $f_S$ is of type IIIa.
In this case, it follows from Theorem~\ref{thm:exclasstypIIIb} that $\pi(r(K))^0$ is a torus,
since $\mathcal{H}^0$ is, so we have $T=K^0$ and ${\rm Z}(r)^0$ is a complex torus. 
But then any element of ${\rm Z}(r)$, hence of $r(\Gamma)$, is semi-simple, in contradiction 
with the second assertion in assumption (ii).
\end{pf}

\begin{example} 
Assume $k=\overline{\Q_\ell}$ is an algebraic closure of 
the field $\Q_\ell$ of $\ell$-adic numbers, $\Gamma$ is the absolute Galois group of $\Q$, 
and $r$ is a geometric ${\rm GSpin}_n(\overline{\Q_\ell})$-valued representation 
in the sense of Fontaine and Mazur.  
The first part of assumption (ii) holds for instance if 
the multiplicity of the Hodge-Tate weight $0$ of $\pi \circ r$ is $<3$. 
Under the assumptions of Proposition~\ref{prop:gspinkmor}, 
it follows that the collection of conjugacy classes 
in ${\rm GSpin}_n(\overline{\Q_\ell})$ of the semi-simplified Frobenius elements of $r$ 
determine $r$ up to conjugacy. In the case $n=7$, this improves~\cite[Prop. 5.2]{kretshin}.
\end{example}

\appendix

\section{}

We gather in this appendix a few lemmas that we used.
The first is the following folklore variant of the acceptability of ${\rm O}(n)$ and ${\rm U}(n)$. 

\begin{lemma}
\label{lem:straccclassical}
Let $G$ be either ${\rm U}(n)$ or ${\rm O}(n)$ for $n\geq 1$, $\rho: G \rightarrow {\rm GL}_n(\C)$ its tautological representation, and $r,r' : \Gamma \rightarrow G$ two group morphisms. Then 
$r$ and $r'$ are conjugate in $G$ if, and only if, the representations $\rho \circ r$ and $\rho \circ r'$ of $\Gamma$ are isomorphic.  Moreover, the same result holds if $G$ is ${\rm SU}(n)$, or ${\rm SO}(n)$ with $n$ odd.
\end{lemma}

\begin{pf} 
Two elements $g,g' \in G$ are conjugate in $G$ if, and only if, $\rho(g)$ and $\rho(g')$ have the same characteristic polynomial. The non trivial implication of the statement is then equivalent to the acceptability of $G$ (proved {\it e.g} in \cite{larsen1}).
\end{pf}

The next lemma is about the transfer morphism to an index $2$ subgroup.

\begin{lemma} 
\label{lem:trancar}
Let $\Gamma$ be a group, $\chi : \Gamma \rightarrow \{ \pm 1\}$ an order $2$ character, $\Gamma_0 \subset \Gamma$ the kernel of $\chi$, $c$ a character of $\Gamma_0$ and $t$ the transfer of $c$ to $\Gamma$. Then:
\begin{itemize}
\item[(i)] For all $\gamma \in \Gamma_0$ and
$z \in \Gamma \smallsetminus \Gamma_0$, we have  $t(z) = c(z^2)$ and 
$t(\gamma)= c(\gamma z^{-1} \gamma z)$. \ps
\item[(ii)] If $U$ is a finite dimensional representation of 
$\Gamma_0$ with determinant $c$, then $\det {\rm Ind}_{\Gamma_0}^\Gamma U = \chi^{\dim U} t$. 
\end{itemize}
\end{lemma}

\begin{pf} 
Part (i) is straightforward and part (ii) is due to Gallagher \cite{gallagher}.
\end{pf}

 The second is about a notion of {\it orthogonal induction}.

\begin{lemma} 
\label{lem:realind}
Let $V$ be an Euclidean space, $\Gamma$ a group, $\rho : \Gamma \rightarrow {\rm O}(V)$ a representation, $\Gamma_0 \subset \Gamma$ an index $2$ subgroup and $z \in \Gamma \smallsetminus \Gamma_0$. Assume that there is a $\Gamma_0$-stable subspace $V_0 \subset V$ such that:\begin{itemize}
\item[(i)] $V_0$ is a direct sum of absolutely irreducible representations of $\Gamma_0$.\ps
\item[(ii)] $V = V_0 \oplus zV_0$, or equivalently, the natural morphism ${\rm Ind}_{\Gamma_0}^\Gamma V_0 \rightarrow V$ is an isomorphism.\ps
\end{itemize}
Then there is a $\Gamma_0$-stable subspace $U_0 \subset V$ which is isomorphic to $V_0$
as $\Gamma_0$-module, and satisfying $V = U_0 \perp zU_0$.
\end{lemma}

\begin{pf}
Consider first the case where $V_0$ is an absolutely irreducible $\R[\Gamma_0]$-module.
We have the $\Gamma_0$-stable decompositions $V = V_0 \oplus z V_0$ and $V = V_0 \perp V_0^\perp$. If $zV_0$ is not isomorphic to $V_0$, the orthogonal projection $z V_0 \rightarrow V_0$
is zero, so we have $zV_0 = V_0^\perp$ and we are done. Otherwise, the $\R[\Gamma_0]$-module $V$ is isotypical. As we have ${\rm End}_{\R[\Gamma_0]}(V_0)=\R$ by assumption, and by the theory of isotypic components, we may assume that $V$ is the tensor product of $V_0$ and of some Euclidean plane $P \simeq \R^2$, and that the action of $\Gamma_0$ on $V= V_0 \otimes P$ is the given one on the first factor, and trivial on the second. \ps

The centralizer of $\rho(\Gamma_0)$ in ${\rm O}(V)$ is $1 \otimes {\rm O}(P)$, and that of $1 \otimes {\rm O}(P)$ is ${\rm O}(V_0) \otimes 1$. The element $\rho(z)$ acts on $V$ by normalizing $\rho(\Gamma_0)$, hence by normalizing $1 \otimes {\rm O}(P)$ as well. As each automorphism of ${\rm O}(P)$ is inner, we may thus write $\rho(z) = \gamma \otimes \delta$ for some $\gamma \in {\rm O}(V_0)$ and $\delta \in {\rm O}(P)$.
As $z^2  \in \Gamma_0$, we have $$\delta^2  \in ({\rm O}(V_0) \otimes 1)  \cap (1 \otimes {\rm O}(P)) = \{\pm {\rm id}_V\}.$$  The proper $\Gamma_0$-stable subspaces of $V$ are the $V_0 \otimes v$ for $v \in P$ nonzero. By assumption (ii), there is $v \in P$
such that $\delta(v) \notin \R v$, {\it i.e.} $\delta$ is not a homothety. It follows that either $\delta$ is an orthogonal symmetry (case $\delta^2=1$), or a rotation of angle $\pi/2$ (case $\delta^2=-1$). In both cases there is a nonzero $v_0 \in P$ such that $v_0$ and $\delta(v_0)$ are orthogonal. The $\R[\Gamma_0]$-module $U_0 = V_0 \otimes v_0$ does the trick.\ps

Consider now the general case. 
Let $A$ be an irreducible $\R[\Gamma_0]$-submodule of $V_0$. 
We have $zA \cap A =\{0\}$ by (ii).
By applying the first paragraph to $V_1=A \oplus zA$, we may find a $\Gamma_0$-stable 
$A' \subset V_1$ isomorphic to $A$ and with $V_1 = A' \perp z A'$. 
Write $V=V_1 \perp V_2$; both $V_i$ are $\Gamma$-stable. 
Let $B$ be an $\R[\Gamma_0]$-module such that
$V_0 \simeq A \oplus B$. 
We must have $V_2 \simeq {\rm Ind}_{\Gamma_0}^\Gamma B$ by semi-simplicity.
By induction on $\dim V$, we may write $V_2 = B' \perp zB'$ with $B'$ a $\Gamma_0$-stable subspace of $V_2$ isomorphic to $B$ as $\Gamma_0$-module. 
The subspace $U_0 = A' \perp B'$ concludes the proof. \end{pf}

We also used the more specific:

\begin{lemma} 
\label{lem:propF}
Assume we are in the situation of Proposition~\ref{prop:devtypIII}. Then:
\begin{itemize}
\item[(i)] the $\R[\Gamma_0]$-module $F$ does not contain $1$ nor $\eta_{|\Gamma_0}$,\ps
\item[(ii)] the $\R[\Gamma_0]$-module $V_0$ is a direct sum of absolutely irreducible representations.\end{itemize}
\end{lemma} 

\begin{pf} 
Assume that the trivial representation $1_0$ of $\Gamma_0$ appears in $F$.
Recall that the trivial representation $1$ of $\Gamma$ does not appear in $F$ by definition in types II or III. If $1_0$ appears in $V_0$ (or equivalently, in its outer conjugate by $\Gamma/\Gamma_0$), 
then  ${\rm Ind}_{\Gamma_0}^\Gamma 1_0 \simeq 1 \oplus \chi$ embeds in $F$, a contradiction.
So we are in type IIIa and $1_0$ appears in $Q_{|\Gamma_0}$. But $\det Q$ is $1$ on $\Gamma_0$,
so we have $Q_{|\Gamma_0} \simeq 1_0 \oplus 1_0$ and again $Q \simeq 1 \oplus \chi$. 
For similar reasons, $\eta_{|\Gamma_0}$ does not occur in $F$:  we have $$ {\rm Ind}_{\Gamma_0}^\Gamma \eta_{|\Gamma_0} \simeq \eta \otimes {\rm Ind}_{\Gamma_0}^\Gamma 1 \simeq \eta \oplus \eta \chi$$
and neither $\eta$ nor $\eta \chi$ appears in $F$ as $(r,\eta)$ is unacceptable. This proves (i). 
If $S$ is an irreducible $\R[\Gamma_0]$-submodule of $V_0$, we have $1 \leq \dim S \leq 3$.
If $S$ is not absolutely irreducible, we necessarily have $\dim S=2$ and ${\rm End}_{\R[\Gamma_0]} S = \C$. But in this case we have $\det S = 1$. As $\det V_0 = \eta_{|\Gamma_0} \neq 1$, we have $\dim V_0=3$ and so $V_0 \simeq S \oplus \eta_{|\Gamma_0} $, in contradiction with (i).
\end{pf}


\begin{thebibliography}{99}

\bibitem[\textsc{Ada}96]{adams} J. F. Adams, \textit{Lectures on exceptional Lie groups}, Chicago Lectures of Mathematics (1996).\ps\ps

\bibitem[\textsc{Bla}94]{blasius} D. Blasius, {\it On multiplicities for ${\rm SL}(n)$}, Israel Journal of Math. 88, 237--251 (1994).\ps\ps

\bibitem[\textsc{BtD}]{tdieck} T. Br\"ocker \& T. tom Dieck, {\it Representations of Compact Lie Groups}, Springer GTM 98 (1985).\ps\ps 

\bibitem[\textsc{Che}19]{chg2} G. Chenevier, {\it Subgroups of ${\rm Spin}(7)$ or ${\rm SO}(7)$ 
with each element conjugate to some element of ${\rm G}_2$ 
and applications to automorphic forms}, Doc. Math. 24, 95--161 (2019).\ps\ps

\bibitem[\textsc{CG}18]{letter} G. Chenevier \& W. T. Gan, {\it Letter to Larsen} (2018).\ps\ps

\bibitem[\textsc{Gal}65]{gallagher} P. X . Gallagher, {\it Determinants of representations of finite groups},
Abh. math. Seminar Univ. Hamburg  3/4,  162--167 (1965).\ps\ps

\bibitem[\textsc{GS}22]{gansavin} W. T. Gan \& G. Savin, {\it The Local Langlands Conjecture for ${\rm G}_2$},
\href{https://arxiv.org/abs/2209.07346}{arXiv} preprint (2022).\ps\ps

\bibitem[\textsc{Gri}95]{griess} R. Griess, \textit{Basic conjugacy theorems for ${\rm G}_2$}, Inventiones Math. 121, 257--278 (1995). \ps\ps

\bibitem[\textsc{KS}]{kretshin} A. Kret \& S. W. Shin, {\it Galois representations for general symplectic groups}, to appear in J. Eur. Math. Soc.  \ps \ps

\bibitem[\textsc{Hum}98]{humphreys} J. E. Humphreys, {\it Linear Algebraic Groups}, Springer Verlag, 
Graduate Textes in Mathematics 21 (1998).\ps\ps

\bibitem[\textsc{Lar}94]{larsen1} M. Larsen, {\it On the conjugacy of element-conjugate homomorphisms}, Israel Journal of Math. 88, 253-277 (1994). \ps\ps

\bibitem[\textsc{Lar}96]{larsen2} M. Larsen, {\it On the conjugacy of element-conjugate homomorphisms II},  Quart. J. Math. Oxford Ser. (2) 47 (1996), no. 185, 73--85. \ps\ps

\bibitem[\textsc{Tat79}]{tate} J. Tate, \textit{Number theoretic background}, Proc. of Symposia in Pure Math. 33 vol II, Corvallis Conference, Oregon 1977, American Math. Soc., 3--26 (1979).\ps\ps

\bibitem[\textsc{Wan}15]{wang} S. Wang, {\it On local and global conjugacy}, Journal of Algebra 439, 334-359 (2015).\ps\ps

\bibitem[\textsc{Wei}74]{weil} A. Weil, {\it Exercices dyadiques}, Inventiones Mathematicae 27, 1--22 (1974).\ps\ps

\bibitem[\textsc{Yu}21]{yu} J. Yu, {\it Acceptable compact Lie groups}, Peking Mathematical Journal 5, 427--446 (2022).\ps\ps
\end{thebibliography}
\end{document}